\documentclass{article}

\usepackage{cmd}
\usepackage{format}

\title{Iterating Semi-proper Forcing using Virtual Models}
\author{Obrad Kasum\footnote{Mr. Obrad Kasum is co‑funded by the European Commission under the programme H2020‑MSCA‑COFUND‑2019 Grant agreement 945332 \includegraphics[height=2.5mm]{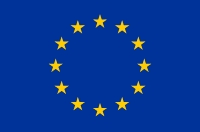}}, Boban Veli\v ckovi\' c}
\date{March 22, 2023}

\begin{document}

\maketitle
\begin{abstract}
    By a virtual model, we mean a model of set theory which is elementary in its transitive closure.
    Virtual models are first used by Neeman \cite{neeman2014forcing} to iterate forcing.
    That paper is concerned with proper forcing.
    The method was then adjusted by Veličković to iterate the case of semi-proper forcing and this was drafted in \cite{velickovic2021iteration}.
    We here straighten the details and futher elaborate on Veličković's method.
    The first section collects facts about virtual model, the second section describes the iteration, and the third one illustrates the method in the case of getting saturation of $\ns$ (loosely relying on \cite{schindler2016nsomega1}).
\end{abstract}
\tableofcontents
\newpage

\section{Virtual Models}
\subsection{Admissible Structures}
\begin{enumerate}
\item\summary We consider a basic notion of an admissible structure.
It will depend on parameter $\mathcal{L}$, where $\mathcal{L}$ is a recursively enumerable first order language containing $\in$.
The language $\mathcal{L}$ will most often be kept implicit.

\item\definition An \textit{admissible structure} $\A:=(A,\in,\dots)$ in language $\mathcal{L}$ is a transitive structure which satisfies $\mathsf{ZFC}^-$ in the extended language.

\item\remark We will often say ``$A$ is admissible" when we really mean ``$\A$ is admissible".

\item\example If $A$ is transitive model for $\mathsf{ZFC}^-$, then $(A,\in)$ is admissible.

\item\example If $A$ is admissible and if $U\subseteq A$ is definable with parameters over $(A,\in)$, then $(A,\in,U)$ is admissible.

\item\example If $\theta>\omega$ is regular, then every structure $(H_\theta,\in,\dots)$ is admissible.
\end{enumerate}

\subsection{Virtual Models: Definition}
\begin{enumerate}
\item\summary We introduce now a generalization of an admissible structure which we call a virtual model.
This notion will play a crucial role in the main construction.

\item\definition A \textit{virtual model (in language $\mathcal{L}$)} is a structure $(M,\in,\dots)$ such that there exists an admissible $\A$ with
$$(M,\in,\dots)\prec\A.$$
    
\item\example Every admissible structure is also a virtual model.

\item\remark If $M\prec\A$, where $\A$ is admissible, then the structure on $M$ is uniquely determined by structure $\A$.

\item\goal We will now show that given virtual model $M$, there exists unique minimal structure $\widehat{M}$ witnessing the fact that $M$ is a virtual model.

\item\lemma Let $ A$ be admissible and let $B\subseteq A$ be transitive.
Then $B\prec_0 A$.

\item\lemma Let $ A$ be admissible and let $M\prec_0 A$ be such that $\trcl(M)=A$.
Then $M\prec A$.

\item\remark We treat function symbols as relation symbols.
In particular, if $ A$ is admissible and $M\subseteq A$ contains all the constants, then $M$ inherits a structure from $A$.
In that case, the inherited structure is the default structure on $M$.

\item\lemma Let $ A$ be admissible and let $M\prec A$.
Then $\trcl(M)$ is admissible and 
$$M\prec\trcl(M)\prec A.$$

\item\proposition Let $M$ be a virtual model.
Then there exists unique structure $\widehat{M}$ on $\trcl(M)$ satisfying $M\prec\widehat{M}$.
Furthermore, for every admissible $A$, it holds that $M\prec A$ if and only if $\widehat{M}\prec A$.
\begin{proof}
Existence follows from the previous lemma.
To verify uniqueness, consider a relation $R$ of structure $M$ and note that
$$R^{\widehat{M}}=\bigcup_{\xi<\ord\cap\widehat{M}}(R^{\widehat{M}}\cap (\widehat{M}\upharpoonright\xi))=\bigcup_{\xi<\ord\cap{M}}(R^{\widehat{M}}\cap (\widehat{M}\upharpoonright\xi))=\bigcup_{\xi<\ord\cap{M}}(R\cap V_\xi)^M$$
which depends only on $M$.
\end{proof}
\end{enumerate}

\subsection{Elementary Rank-initial Segments}
\begin{enumerate}
\item\definition Let $A$ be admissible.
Then $\e_A:=\{\alpha<\ord^A : A\upharpoonright\alpha\prec A\}$.

\item\remark Set $\e_A$ is closed in $\ord^A$.
If $A=V_\kappa$ where $\kappa$ is inaccessible, then $\e_A$ is club in $\kappa$.

\item\proposition Let $A,B$ be admissible and let $\alpha\in\e_A\cap\e_B$.
Suppose that $A\upharpoonright\alpha=B\upharpoonright\alpha$.
Then 
$$\e_A\cap [0,\alpha]=\e_B\cap [0,\alpha]\in A\cap B.$$

\item\definition Let $M$ be a virtual model.
Then $\e_M:=\{\alpha\in\ord\cap M : M\cap V_\alpha\prec M\}.$

\item\remark Every admissible structure is also a virtual model.
In that case, the two definitions coincide.

\item\proposition Let $M$ be a virtual model.
Then $\e_M=\e_{\widehat{M}}\cap M$.

\item\corollary Let $M$ be a virtual model and let $\gamma$ be a limit point of $\e_M$.
Then $\gamma\in\e_{\widehat{M}}$.
\begin{proof}
Follows from the fact that $\e_M\subseteq\e_{\widehat{M}}$ and the fact that $\e_{\widehat{M}}$ is closed.
\end{proof}

\item\label{722}\lemma Let $M$ be a virtual model and let $\alpha\in\e_{\widehat{M}}$.
Suppose that $\beta:=\min(M\cap [\alpha,\infty))$.
Then $\beta\in\e_M$.
\begin{proof0}
    \item We apply Tarski-Vaught test.
    Suppose $a\in M\cap V_\beta$ and $M\models\exists y\phi(a,y)$.
    We want to find $b\in M\cap V_\beta$ such that $M\models\phi [a,b]$.
    
    \item By elementarity, $\widehat{M}\models\exists y\phi(a,y)$.
    We also have $a\in M\cap V_\beta\subseteq M\cap V_\alpha\subseteq\widehat{M}\upharpoonright\alpha$.
    
    \item Since $\alpha\in\e_{\widehat{M}}$, there exists $y\in \widehat{M}\upharpoonright\alpha$ such that $\widehat{M}\models\phi [a,y]$.
    
    \item Since $\alpha\leq \beta$, we have
    $$\widehat{M}\models(\exists y\in V_{\beta})\phi(a,y).$$
    By elementarity,
    $$M\models(\exists y\in V_{\beta})\phi(a,y).$$
    
    \item Thus, there is $b\in M\cap V_\beta$ such that $M\models\phi [a,b]$.
\end{proof0}

\item\label{724}\label{eM lemma}\proposition Let $M$ be a virtual model.
Suppose that $\gamma$ is a limit point of $\e_{\widehat{M}}$ and a limit point of $M\cap\ord$.
Then $\gamma$ is a limit point of $\e_M$.
\begin{proof0}
    \item Let $\alpha<\gamma$ be arbitrary.
    We want to find $\beta\in\e_M\cap\gamma$ such that $\alpha<\beta$.
    
    \item There is $\alpha'\in\e_{\widehat{M}}\cap\gamma$ such that $\alpha<\alpha'$.
    
    \item\label{753} There is $\alpha''\in M\cap\gamma$ such that $\alpha'<\alpha''$.
    
    \item By \ref{722}, we have $\beta:=\min(M\cap [\alpha',\infty))\in\e_{M}$.
    
    \item By \ref{753}, we have $\beta\in [\alpha',\alpha'']\subseteq (\alpha,\gamma)$.
\end{proof0}
\end{enumerate}

\subsection{Reduction of a Virtual Model}
\begin{enumerate}
\item\summary We introduce a basic construction of reducing the amount of information captured by a virtual model.

\item\definition Let $M$ be a virtual model and let $X$ be an arbitrary set.
We define
$$\hull(M,X):=\{f(x) : f:d\rightarrow r,\, f\in M,\, x\in X^{<\omega}\cap d\}.$$

\item\proposition $M$ be a virtual model, let $\alpha:=\sup(\ord\cap M)$, and let $X$ be an arbitrary set.
Then the following holds.
\begin{parts}
    \item $M\prec\hull(M,X)\prec\widehat{M}$
    
    \item $\widehat{\hull(M,X)}=\widehat{M}$
    
    \item $X\cap\widehat{M}\subseteq\hull(M,X)$
    
    \item $\hull(M,X)=\hull(M,X\cap\widehat{M})$
    
    \item For every virtual model $N$ satisfying $M\prec N$ and $X\cap\widehat{M}\subseteq N$, we have $\hull(M,X)\prec N$.
    
    \item $|\hull(M,X)|\leq |M|+|X|$
\end{parts}\qed

\item\definition \textit{The $\alpha$-reduction $M\downarrow\alpha$} of virtual model $M$ is defined as follows: let $\pi:\hull(M,V_\alpha)\rightarrow N$ be the transitive collapse and set $M\downarrow\alpha:=\pi [M]$.

\item\remark We do not require $V_\alpha\subseteq \widehat{M}$ or even $\delta:=\sup(\ord\cap M)\geq\alpha$.
What we do is simply collapsing $M$ while ``freezing" its part below rank $\alpha$.
In particular, if $\delta\leq\alpha$, this means that $M\subseteq V_\alpha$ and $\hull(M,V_\alpha)=\widehat{M}$, i.e. $M\downarrow\alpha= M$.

\item\exercise Let $M$ be a virtual model and let $\alpha$ be an ordinal.
Then $\widehat{M\downarrow\alpha}$ is the transitive collapse of $\hull(M,V_\alpha)$.

\item\definition A virtual model $M$ is said to be \textit{$\alpha$-generated} if $\widehat{M}=\hull(M,V_\alpha)$.

\item\label{col and sigma0}\lemma Let $N$ be a virtual model and let $\pi:N\rightarrow\overline{N}$ be the transitive collapse.
Then for every $\Sigma_0$ formula $\phi(\overline{x})$, we have
$$(\forall\overline{x}\in N^{<\omega})(\phi(\overline{x})\iff\phi(\pi(\overline{x}))).$$
\begin{proof}
\begin{eqnarray}
\phi(\overline{x}) &\iff& \widehat{N}\models\phi(\overline{x}) \label{204}\\
&\iff& N\models\phi(\overline{x}) \label{205}\\
&\iff& \overline{N}\models\phi(\pi(\overline{x}))\label{206}\\
&\iff& \phi(\pi(\overline{x})),\label{207}
\end{eqnarray}
where:
\begin{itemize}
    \item[(\ref{204})] follows by $\Sigma_0$-absoluteness;
    
    \item[(\ref{205})] follows by elementarity;
    
    \item[(\ref{206})] follows since $\pi$ is an isomorphism;
    
    \item[(\ref{207})] follows by $\Sigma_0$-absoluteness.
\end{itemize}
\end{proof}

\item\label{146}\proposition Let $M$ be a virtual model.
Then $M\downarrow\alpha$ is $\alpha$-generated.

\item\label{reductions}\proposition Let $M$ be a virtual model and let $\beta\leq\alpha$.
Then $(M\downarrow\alpha)\downarrow\beta=M\downarrow\beta$.
\begin{proof}
This is a straight forward computation; use \ref{col and sigma0} to show that collapses agree with the computations of $\hull$'s.
\end{proof}

\end{enumerate}

\subsection{\texorpdfstring{$\alpha$}{alpha}-isomorphism}
\begin{enumerate}
\item\summary The relation introduced in this part has for the goal to capture the idea that two virtual models carry the same information up to $\alpha$.

\item\definition Let $M,N$ be virtual models.
We define $M\cong_\alpha N$ to hold if there is an isomorphism $f:\hull(M,V_\alpha)\cong\hull(N,V_\alpha)$ satisfying $f[M]=N$.

\item\proposition Let $M$ be a virtual model.
Then $M\cong_\alpha M\downarrow\alpha$.

\item\proposition Let $M,N$ be $\alpha$-generated virtual models.
If $M\cong_\alpha N$, then $M=N$.

\item\proposition Let $M,N$ be virtual models.
Then $M\cong_\alpha N$ if and only if $M\downarrow\alpha=N\downarrow\alpha$.

\item\corollary $\cong_\alpha$ is an equivalence relation between virtual models.
\end{enumerate}

\subsection{Comparing Virtual Models}
\begin{enumerate}
\item\summary

\item\definition Let $M,N$ be countable virtual models.
Then we define $M\lhd_\alpha N$ to hold if there exists $M'\in N$ such that $M'$ with the inherited structure is a virtual model and $M\cong_\alpha M'$.

\item\definition A virtual model $M$ is said to be \textit{$\xi$-strong} if $V_\xi\subseteq\widehat{M}$.

\item\label{remaining countable}\lemma Let $M$ be an $(\omega+1)$-strong virtual model and let $N\in M$ be a countable virtual model.
Then $|N|^M=\omega$.
\begin{proof0}
    \item Let $\pi:M\rightarrow\overline{M}$ be the transitive collapse of $(M,\in)$.
    We have $\overline{M}\in H_{\omega_1}$ and there is a surjection $f:\omega\rightarrow Y$.
    
    \item We have $f\in H_{\omega_1}\subseteq\widehat{N}$.
    
    \item Thus, $\widehat{N}\models(\overline{M}$ is countable$)$.
    
    \item By absoluteness of transitive collapse, $\widehat{N}\models (\overline{M},\in)\cong (M,\in)$.
    
    \item Thus, $\widehat{N}\models |M|=|\overline{M}|=\omega$ and consequently $N\models|M|=\omega$.
\end{proof0}

\item\proposition Suppose that $\alpha>\omega$.
Then relation $\lhd_\alpha$ is a partial order between countable $(\omega+1)$-strong virtual models.
\begin{proof0}
    \item Suppose that $M,N,P$ are countable $(\omega+1)$-strong virtual models satisfying $M\lhd_\alpha N\lhd_\alpha P$.
    
    \item There exist a virtual model $M'\in N$, an isomorphism 
    $$f:\hull(M,V_\alpha)\cong\hull(M',V_\alpha)$$
    with $f[M]=M'$, a virtual model $N'\in P$, and an isomorphism 
    $$g:\hull(N,V_\alpha)\cong\hull(N',V_\alpha)$$
    with $g[N]=N'$.
    
    
    \item By \ref{remaining countable}, we have $M'\subseteq N$ and $g(M')=g[M']$.
    
    \item For $h\in N'$ and $x\in\dom(f)\cap V_\alpha^{<\omega}$, we have 
    $$g(h(x))=g(h)(g(x))=g(h)(x).$$
    We conclude 
    $$g[\hull(M',V_\alpha)]=\hull(g[M'],V_\alpha)=\hull(g(M'),V_\alpha).$$
    
    \item Thus, $g(M')\in P$ and
    $$g\circ f:\hull(M,V_\alpha)\cong \hull(g(M'),V_\alpha).$$
    
    \item We also verify
    $$(g\circ f)[M]=g[f[M]]=g[M']=g(M'),$$
    which yields the conclusion $M\lhd_\alpha P$.
\end{proof0}

\item\label{alpha-in and downarrow}\proposition Suppose that $M,M',N,N'$ are countable virtual models such that $M\cong_\alpha M'$ and $N\cong_\alpha N'$.
Then $M\lhd_\alpha N$ if and only if $M'\lhd_\alpha N'$.\qed 

\item\label{in gives alpha-in}\proposition Suppose that $M,N$ are countable virtual models and $M\in N$.
Then $M\lhd_\alpha N$.\qed 

\item\label{alpha-in to in}\proposition Let $\alpha$ be an ordinal and let $M,N$ be countable virtual models.
Suppose that $M,\alpha\in N$, that $M$ is $\alpha$-generated, and $M\lhd_\alpha N$.
Then $M\in N$.
\begin{proof0}
    \item There is a virtual model $M'\in N$ such that $M\cong_\alpha M'$.
    
    \item Since $M$ is $\alpha$-generated, we have $M=M'\downarrow\alpha\in N$.
\end{proof0}

\item\label{alpha-in to in for H}\proposition Let $\alpha$ be a beth-fixed-point and let $M,N$ be countable $\alpha$-strong virtual models such that $M\lhd_\alpha N$.
Let $\theta\in M\cap V_\alpha$ be regular and uncountable.
Then $M\cap H_\theta\in N\cap H_\theta$.
\begin{proof}
Since $N$ is $\alpha$-strong and $\theta<\alpha$, we have that $H_\theta=H_\theta^N\in N$.
Let $M'\in N$ be such that $M\cong_\alpha M'$.
We have $M\cap H_\theta=M'\cap H_\theta\in N$.
On the other hand, since $M\cap H_\theta\in [H_\theta]^\omega$, we conclude $M\cap H_\theta\in H_\theta$.
\end{proof}
\end{enumerate}

\subsection{Forcing Extensions of Virtual Models}
\label{Forcing Extensions of Virtual Models}
\begin{enumerate}
    \item\definition Let $A$ be admissible, let $\P\in A$ be a poset, and let $M\prec A$ with $\P\in M$.
    For $G\leadsto A^\P$, we define
    $$M[G]:=\{\tau^G : \tau\in M^\P=A^\P\cap M\},\quad M^G:=M[G]\cap A.$$
    
    \item\label{vm and genericity}\proposition  Let $A$ be admissible, let $\P\in A$ be a poset, let $M\prec A$ with $\P\in M$, and let $G\leadsto A^\P$.
    Then the following holds.
    \begin{parts}
        \item $(M[G],\in,M^G,\P,G)\prec (A[G],\in,A,\P,G)$
        
        \item $\widehat{M[G]}=\widehat{M}[G]$
        
        \item $M\prec M^G\prec A$
        
        \item\label{329} $M^G[G]=M[G]$
        
        \item\label{338} $(M^G)^G=M^G$.
    \end{parts}
    \qed
    
    \item\label{strength and forcing}\proposition Let $M$ be an $\alpha$-strong virtual model, let $\alpha$ be such that $\widehat{M}\upharpoonright\alpha\models\mathsf{ZFC}^-$, let $\P\in M\cap V_\alpha$, and let $G\leadsto V^\P$.
    Then $(\widehat{M}\upharpoonright\alpha)[G]\subseteq\widehat{M[G]}$.\qed
    
    \item\label{hulls and gen}\lemma Let $M$ be a virtual model, let $\alpha\in\widehat{M}$ satisfy $\widehat{M}\upharpoonright\alpha\models\mathsf{ZFC}^-$, let $\P\in M\cap V_\alpha$, and let $G\in V$ satisfy $G\leadsto \widehat{M}^\P$.
    Then $\hull(M[G],V_\alpha)=\hull(M,V_\alpha)[G]$.
    \begin{proof0}
        \item We first consider inclusion ($\supseteq$).
        \begin{pfenum}
        
        \item Let $\dot{y}\in\hull(M,V_\alpha)$ and let us verify $\dot{y}^G\in\hull(M[G],V_\alpha)$.
        
        \item We have $\dot{y}=f(x)$ for a function $f\in M$ and an $x\in V_\alpha\cap\dom(f)$.
        We may choose $f$ so that $\im(f)\subseteq\widehat{M}^\P$.
        
        \item Let $g\in M[G]$ be the function satisfying $\dom(g)=\dom(f)$ and $g(u):=f(u)^G$.
        We have $g\in M[G]$ and $x\in\dom(g)\cap V_\alpha$.
        
        \item Hence, $\dot{y}^G=f(x)^G=g(x)\in\hull(M[G],V_\alpha)$.
        \end{pfenum}
        
        \item We consider  now inclusion ($\subseteq$)
        \begin{pfenum}
        
        \item Let us consider $x_0\in\hull(M[G],V_\alpha)$ and let us show $x_0\in\hull(M,V_\alpha)[G]$.
        
        \item There exist a function $f\in M[G]$ and $x\in V_\alpha\cap\dom(f)$ such that $x_0=f(x)$.
        
        \item There exist $\P$-names $\dot{f},\dot{d}$ such that 
        $$\widehat{M}^\P\models(\dot{f}\mbox{ is a function on }\dot{d})$$
        and $\dot{f}^G=f$.
        
        \item Note that
        $$x\in\widehat{M[G]}\upharpoonright\alpha=\widehat{M}[G]\upharpoonright\alpha=(\widehat{M}\upharpoonright\alpha)[G],$$
        so there exists a $\P$-name $\dot{x}\in \widehat{M}\upharpoonright\alpha$ such that $\dot{x}^G=x$.
        
        \item Let $\xi\in M$ satisfy $\xi>\alpha$ and let $e:=\{\dot{y}\in V_\xi^M : \dot{y}$ is a $\P$-name$\}\in M$.
        
        \item There exists a function $g$ such that $\dom(g)=e$ and
        $$(\forall\dot{y}\in e)(\forall p\in\P)(p\Vdash_\P^{\widehat{M}}\dot{y}\in\dot{d}\implies p\Vdash_\P^{\widehat{M}}\dot{f}(\dot{y})=\ulcorner g(\dot{y})\urcorner).$$
        By elementarity, function $g$ can be chosen in $M$.
        
        \item We have $\dot{z}:=g(\dot{x})\in\hull(M,V_\alpha)$ and $\dot{z}^G\in\hull(M,V_\alpha)[G]$.
        
        \item Let $p\in G$ be such that $p\Vdash_\P^{\widehat{M}}\dot{x}\in\dot{d}$.
        Since $\dot{x}\in e$, we have $p\Vdash_\P^{\widehat{M}}\dot{f}(\dot{x})=\dot{z}$ and consequently $\dot{f}^G(\dot{x}^G)=\dot{z}^G$.
        
        \item Hence, $x_0=f(x)=\dot{f}^G(\dot{x}^G)=\dot{z}^G\in\hull(M,V_\alpha)[G].$
        \end{pfenum}
    \end{proof0}
    
    \item\label{cong-alpha and gen}\proposition Let $M,N$ be virtual models and let $\alpha\in\widehat{M}\cap\widehat{N}$.
    Suppose that 
    $$\widehat{M}\upharpoonright\alpha=\widehat{N}\upharpoonright\alpha\models\mathsf{ZFC}^-$$
    and that $M\cong_\alpha N$.
    Let $\P\in M\cap V_\alpha$ and let $G\in V$ satisfy $G\leadsto\widehat{M}^\P$.
    Then 
    $$\P\in N\cap V_\alpha,\quad G\leadsto\widehat{N}^\P,\quad M[G]\cong_\alpha N[G],\quad M^G\cong_\alpha N^G.$$
    Furthermore, the isomorphism witnessing $M[G]\cong_\alpha N[G]$ extends the isomorphism witnessing $M^G\cong_\alpha N^G$, which in turn extends the isomorphism witnessing $M\cong_\alpha N$.
    \begin{proof0}
        \item $M\cong_\alpha N$ implies $M\cap V_\alpha=N\cap V_\alpha$, so $\P\in N\cap V_\alpha$.
        The fact that $G\leadsto \widehat{N}^\P$ is obvious.
        
        \item We verify below that $M[G]\cong_\alpha N[G]$ and that the isomorphism witnessing this fact extends the isomorphism witnessing $M\cong_\alpha N$.
        Since $M^G$ and $N^G$ are computed as the appropriate grounds, we will immediately have $M^G\cong_\alpha N^G$, via the restricted isomorphism.
        
        \item Let $F:\hull(M,V_\alpha)\cong\hull(N,V_\alpha)$ be such that $F[M]=N$.
        
        \item By \ref{hulls and gen}, every element of $\hull(M[G],V_\alpha)$ is of the form $\dot{x}^G$ for some $\dot{x}\in\hull(M,V_\alpha)^\P$.
        Also, $\dot{y}^G\in \hull(N[G],V_\alpha)$ for every $\dot{y}\in\hull(N,V_\alpha)^\P$.
        
        \item Hence, we can define $F[G]:\hull(M[G],V_\alpha)\rightarrow\hull(N[G],V_\alpha)$ by setting $F[G](\dot{x}^G):=F(\dot{x})^G$ for $\dot{x}\in\hull(M,V_\alpha)^\P$.
        We want to show that $F[G]$ is an isomorphism, that $(F[G])[M[G]]=N[G]$, and that $F[G]\upharpoonright\hull (M,V_\alpha)=F$.
        
        \item\claim For $*\in\{\in,=\}$ and $\dot{x}_0,\dot{x}_1\in\hull(M,V_\alpha)^\P$, we have
        $$\dot{x}_0^G*\dot{x}_1^G\iff F(\dot{x}_0)^G*F(\dot{x}_1)^G.$$
        \begin{proof}
        \begin{eqnarray}
        \dot{x}_0^G*\dot{x}_1^G & \iff & (\exists p\in G)\hull(M,V_\alpha)\models(p\Vdash\dot{x}_0*\dot{x}_1)\\
        &\iff& (\exists p\in G)\hull(N,V_\alpha)\models(p\Vdash F(\dot{x}_0)*F(\dot{x}_1))\label{364}\\
        &\iff& F(\dot{x}_0)^G*F(\dot{x}_1)^G,
        \end{eqnarray}
        where (\ref{364}) follows from $F\upharpoonright V_\alpha=\id$ and $\P\in V_\alpha$.
        \end{proof}
        
        \item Hence, $F[G]$ is correctly defined injection agreeing with $\in$.
        
        \item Analogously, we have injection 
        $$F^{-1}[G]:\hull(N[G],V_\alpha)\rightarrow\hull(M[G],V_\alpha)$$
        agreeing with $\in$.
        
        \item Clearly, $F[G]\circ F^{-1}[G]=\id$ and $F^{-1}[G]\circ F[G]=\id$, leading to the conclusion that 
        $$F[G]:\hull(M[G],V_\alpha)\cong\hull(N[G],V_\alpha).$$
        
        \item For $x\in\hull(M,V_\alpha)$, we have
        $$F[G](x)=F[G](\check{x}^G)=F(\check{x})^G=((F(x))^{\check{}}\,)^G=F(x).$$
        
        \item Hence, it remains to establish $(F[G])[M[G]]=N[G]$, i.e. $(F[G])[M[G]]\subseteq N[G]$ and $(F^{-1}[G])[N[G]]\subseteq M[G]$.
        
        \item We consider here only the first conjuct.
        Let $y\in M[G]$ be arbitrary.
        
        \item By definition of $M[G]$, we have $y=\dot{x}^G$ for some $x\in M^\P\subseteq\hull(M,V_\alpha)^\P$.
        
        \item By definition of $F[G]$, we have $F[G](y)=F(\dot{x})^G$.
        
        \item Since $\dot{x}\in M$ and $F[M]=N$, we have $F(\dot{x})\in N^\P$ and consequently $F[G](y)=F(\dot{x})^G\in N[G]$.
    \end{proof0}
    
    \item\label{alpha-in and gen ext}\corollary Let $M,N$ be countable virtual models and let $\alpha\in\widehat{M}\cap\widehat{N}$.
    Suppose that 
    $$\widehat{M}\upharpoonright\alpha=\widehat{N}\upharpoonright\alpha\models\mathsf{ZFC}^-$$
    and that $M\lhd_\alpha N$.
    Let $\P\in M\cap V_\alpha$ and let $G\in V$ satisfy $G\leadsto\widehat{M}^\P$.
    Then $M[G]$ and $N[G]$ are countable and $M[G]\lhd_\alpha N[G]$.
    \begin{proof0}
        \item $M,N$ surject onto $M[G],N[G]$, respectively, so $M[G],N[G]$ are countable.
        
        \item Let $M'\in N$ be a virtual model in $V$ satisfying $M\cong_\alpha M'$.
        
        \item By \ref{cong-alpha and gen}, we have $M[G]\cong_\alpha M'[G]\in N[G]$.
        
        \item Thus, $M[G]\lhd_\alpha N[G]$.
    \end{proof0}
    
    \item\lemma Let $M$ be a virtual model and let $\alpha\in\widehat{M}$ be such that $\widehat{M}\upharpoonright\alpha\models\mathsf{ZFC}^-$.
    Suppose that $M$ is $\alpha$-generated.
    Let $\P\in M\cap V_\alpha$ and let $G\in V$ satisfy $G\leadsto\widehat{M}^\P$.
    Then $M[G]$ and $M^G$ are $\alpha$-generated.
    \begin{proof}
    To see that $M[G]$ is $\alpha$-generated, compute as follows:
    $$\widehat{M[G]}=\widehat{M}[G]=\hull(M,V_\alpha)[G]=\hull(M[G],V_\alpha).$$
    To see that $M^G$ is $\alpha$-generated, recall that $M\prec M^G\prec\widehat{M}$.
    \end{proof}
    
    \item\label{reductions and forcing}\proposition Let $M$ be a virtual model and let $\alpha\in\widehat{M}$ be such that $\widehat{M}\upharpoonright\alpha\models\mathsf{ZFC}^-$, let $\P\in M\cap V_\alpha$, and let $G\in V$ satisfy $G\leadsto\widehat{M}^\P$.
    Then $(M\downarrow\alpha)[G]=M[G]\downarrow\alpha$ and $(M\downarrow\alpha)^G=M^G\downarrow\alpha$.
    \begin{proof}
    Since $M\downarrow\alpha\cong_\alpha M$, we have $(M\downarrow\alpha)[G]\cong_\alpha M[G]$.
    Since $M\downarrow\alpha$ is $\alpha$-generated, so is $(M\downarrow\alpha)[G]$ and we conclude 
    $$(M\downarrow\alpha)[G]= M[G]\downarrow\alpha.$$
    The other part is completely analogous.
    \end{proof}
\end{enumerate}

\subsection{Semi-proper Forcing}
\begin{enumerate}
\item\label{semi-gen}\definition Let $\P$ be a poset, let $\theta>2^{2^{|\trcl(\P)|}}$ be regular, and let $M\prec (H_\theta,\in,\P)$ be countable.
\begin{parts}
    \item For $G\leadsto V^\P$, we say that $G$ is \textit{semi-$M^\P$-generic} if $M[G]\cap\omega_1^V=M\cap\omega_1^V$.
    
    \item For $p\in\P$, we say that $p$ is \textit{semi-$M^\P$-generic} if $p\Vdash M[\G]\cap\omega_1^V=M\cap\omega_1^V$.
\end{parts}

\item\label{semi-gen cond and omega1}\lemma Let $\P$ be a poset, let $\theta>2^{2^{|\trcl(\P)|}}$ be regular, and let $M\prec (H_\theta,\in,\P)$ be countable.
If $p\in\P$ is semi-$M^\P$-generic, then $p\Vdash\omega_1=\omega_1^V$.
\begin{proof0}
    \item Let $g\leadsto V^\P$ with $p\in g$.
    
    \item Since $M[g]\prec H_\theta^{V[g]}$, it holds that ``$\omega_1$ is the smallest ordinal $\xi$ with the property that $\sup(M[g]\cap\xi)<\xi$".
    
    \item Since $M[g]\cap\omega_1^V=M\cap\omega_1^V$ and the last set is an ordinal, we see that the solution for $\xi$ is exactly $\xi=\omega_1^V$.
\end{proof0}

\item\definition Let $\P$ be a poset and let $\theta>2^{2^{|\trcl(\P)|}}$ be regular.
Poset $\P$ is \textit{semi-proper} if for every countable $M\prec (H_\theta,\in,\P)$ and every $p\in\P\cap M$ there exists $q\leq p$ which is semi-$M^\P$-generic.

\item\proposition Suppose that $\P$ is semi-proper.
Then:
\begin{parts}
    \item\label{128} $\P\Vdash\omega_1=\omega_1^V$
    
    \item\label{130} for every stationary $S\subseteq\omega_1$, it holds that $\P\Vdash(S$ is stationary$)$.
\end{parts}
\begin{proof0}
    \item Let $\theta$ be sufficiently large regular.
    
    \item By \ref{semi-gen cond and omega1}, it suffices for \ref{128} to show that
    $$D:=\{p\in\P : (\exists M\prec (H_\theta,\in,\P))(|M|=\omega\wedge p\mbox{ is semi-}M^\P\mbox{-generic})\}$$
    is dense in $\P$.
    
    \item Let $p_0\in\P$ be arbitrary and let $M\prec(H_\theta,\in,\P,p_0)$ be countable.
    
    \item Since $\P$ is semi-proper, there is $p\leq p_0$ which is semi-$M^\P$-generic.
    Clearly, $p\in D$.
    
    \item For \ref{130}, suppose that $S$ is stationary in $\omega_1$, that $\Vdash\dot{f}:\omega_1\rightarrow\omega_1$, and let $p\in\P$ be arbitrary.
    We want to find $q\leq p$ and $\xi\in S$ such that $q\Vdash\dot{f}[\xi]\subseteq\xi$.
    
    \item Set of all countable $M\prec (H_\theta,\in,\P,p,\dot{f})$ is club in $[H_\theta]^\omega$.
    
    \item Since $S\subseteq\omega_1$ is stationary in $\omega_1$, there exists a countable $M\prec (H_\theta,\in,\P,p,\dot{f})$ satisfying $\xi:=M\cap\omega_1\in S$.
    
    \item There exists $q\leq p$ which is semi-$M^\P$-generic.
    
    \item For every $g\leadsto M^\P$ with $q\in g$, we have $\xi=M[g]\cap\omega_1$ and $M[g]\models\dot{f}^g:\omega_1\rightarrow\omega_1$.
    We conclude $\dot{f}^g[\xi]\subseteq\xi$.
\end{proof0}

\item\label{semi-gen vm}\definition Let $\alpha$ be a beth-fixed point, let $M$ be a countable $\alpha$-strong virtual model, and let $\P\in V_\alpha\cap M$ be a poset.
\begin{parts}
    \item For $G\leadsto V^\P$, we say that $G$ is \textit{semi-$M^\P$-generic} if $M[G]\cap\omega_1^V=M\cap\omega_1^V$.
    
    \item For $p\in\P$, we say that $p$ is \textit{semi-$M^\P$-generic} if $p\Vdash M[\G]\cap\omega_1^V=M\cap\omega_1^V$.
\end{parts}

\item\proposition Let $\alpha$ be a beth-fixed point, let $M$ be a countable $\alpha$-strong virtual model, let $\P\in V_\alpha\cap M$ be a poset, and let $\theta>2^{2^{|\trcl(\P)|}}$ be a regular cardinal satisfying $\theta<\alpha$.
For every $G\leadsto V^\P$, we have that $G$ is semi-$M^\P$-generic in the sense of \ref{semi-gen vm} if and only if it is semi-$(H_\theta\cap M)^\P$-generic in the sense of \ref{semi-gen}.\qed 

\item\label{cong-alpha and semi-gen}\proposition Let $\gamma$ satisfy $V_\gamma\models\mathsf{ZFC}^-$, let $M,N$ be countable $\gamma$-strong virtual models, let $\P\in M\cap V_\gamma$, and let $G\leadsto V^\P$.
Suppose that $M\cong_\gamma N$.
Then $G$ is semi-$M^\P$-generic if and only if it is semi-$N^\P$-generic.
\begin{proof}
We have $M[G]\cong_\gamma N[G]$ and consequently
$$M[G]\cap\omega_1^V=N[G]\cap\omega_1^V.$$
\end{proof}

\item\label{reductions and semi-genericity}\corollary Let $\gamma$ satisfy $V_\gamma\models\mathsf{ZFC}^-$, let $M$ be a countable $\gamma$-strong virtual model, let $\P\in M\cap V_\gamma$, and let $p\in\P$.
Then $p$ is semi-$M^\P$-generic if and only if it is semi-$(M\downarrow\gamma)^\P$-generic.\qed 
\end{enumerate}

\subsection{Iterated Forcing Extensions of Virtual Models}
\begin{enumerate}
\item\definition A \textit{forcing iteration (with support $E$)} is a family $\Vec{\P}=(\P_\alpha :  \alpha\in E)$ where:
\begin{parts}
    \item $E\subseteq\ord$ and $E\not=\emptyset$;
    
    \item $\P_\alpha$ is a poset for all $\alpha$;
    
    \item $\P_\alpha$ is a complete sub-poset of $\P_\beta$ for all $\alpha\leq\beta$ in $E$.
\end{parts}

\item\proposition Let $\Vec{\P}$ be an iteration with support $E$ and let $\Q:=\cup\{\P_\alpha : \alpha\in E\}$.
Then $\Vec{\P}^\frown\Q$ is an iteration with support $E\cup\{\sup\{\xi+1 : \xi\in E\}\}$.


\item\definition Let $\Vec{\P}$ be a forcing iteration with support $E$, let $\delta\in E$, let $\alpha\in E\cap\delta$, and let $G_\delta\leadsto V^{\P_\delta}$.
Then the \textit{$\alpha$-restriction of the generic $G_\delta$} is defined as $G_\alpha:=G_\delta\cap\P_\alpha$.

\item\remark It holds that $G_\alpha\leadsto V^{\P_\alpha}$.

\item\definition Let $A$ be admissible and let $\Vec{\P}$ be a forcing iteration with support $E$.
Iteration $\Vec{\P}$ is said to be \textit{point-wise definable over $A$} if $\Vec{\P}\subseteq A$ and for every $\alpha\in E$, constant $\P_\alpha$ is definable over $A$ from the parameter $\alpha$.

\item\definition Let $A$ be admissible and let $\Vec{\P}$ be an iteration with support $E$.
Suppose that $\Vec{\P}$ is point-wise definable over $A$.
Let $M\prec A$, let $\delta\in E$, and let $G_\delta\leadsto V^{\P_\delta}$.
Elementary submodels $M^{G_\delta}_{<\alpha}$ and $M^{G_\delta}_\alpha$ of $A$ and elementary submodels $M_\alpha [G_\delta]$ of $A[G_\alpha]$ (inside $V[G_\alpha]$) are defined by recursion on $\alpha\in E\cap [0,\delta]$, as follows.
\begin{parts}
    \item $M^{G_\delta}_{<\alpha}:=M\cup\bigcup_{\xi\in E\cap\alpha} M^{G_\delta}_\xi$
    
    \item\label{47} if $\alpha\in M^{G_\delta}_{<\alpha}$, we define $M_\alpha [G_\delta]:=M^{G_\delta}_{<\alpha} [G_\alpha]$ and $M^{G_\delta}_\alpha:=(M^{G_\delta}_{<\alpha})^{G_\alpha}=M_\alpha [G_\delta]\cap A$;
    
    \item if $\alpha\not\in M^{G_\delta}_{<\alpha}$, we define $M^{G_\delta}_\alpha:=M^{G_\delta}_{<\alpha}$ and leave $M_\alpha [G_\delta]$ undefined.
\end{parts}

\item\proposition The previous definition is correct and for all $\alpha,\beta\in E$, we have:
\begin{parts}
    \item $M^{G_\delta}_{<\alpha}\prec M^{G_\delta}_\alpha\prec A$
    
    \item $M^{G_\delta}_\alpha\prec M^{G_\delta}_{<\beta}$ whenever $\alpha<\beta$;
    
    \item $M_\alpha [G_\delta]\prec A[G_\alpha]$ whenever $\alpha\in M^{G_\delta}_{<\alpha}$.
\end{parts}
\qed
    
    

\item\proposition Let $A$ be admissible, let $\Vec{\P}$ be an iteration with support $E$ that is point-wise definable over $A$, let $M\prec A$, let $\delta_1,\delta_2\in E$ with $\delta_1\leq\delta_2$, let $G_{\delta_2}\leadsto V^{\P_{\delta_2}}$.
Then for all $\alpha\in E\cap [0,\delta_1]$, it holds that:
\begin{parts}
    \item $M^{G_{\delta_1}}_{<\alpha}=M^{G_{\delta_2}}_{<\alpha}$
    
    \item $M^{G_{\delta_1}}_{\alpha}=M^{G_{\delta_2}}_{\alpha}$
    
    \item $M_\alpha [G_{\delta_1}]=M_\alpha [G_{\delta_2}]$ whenever is either of them defined.
\end{parts}
\qed 

\item\proposition Let $A$ be admissible, let $\Vec{\P}$ be an iteration with support $E$ that is point-wise definable over $A$, let $M\prec A$, let $\delta\in E$, let $G_{\delta}\leadsto V^{\P_{\delta}}$, and let $\alpha\in E\cap [0,\delta]$.
Suppose that $\alpha\in M$.
Then
$$M^{G_\delta}_\alpha=M^{G_\alpha},\quad M_\alpha [G_\delta]=M[G_\alpha].$$
\begin{proof0}
    \item It suffices to prove $M_\alpha [G_\delta]=M[G_\alpha]$.
    
    \item We have $M\subseteq M^{G_\delta}_{<\alpha}$, so we conclude
    $$M[G_\alpha]\subseteq M^{G_\delta}_{<\alpha}[G_\alpha]=M_\alpha [G_\delta].$$
    
    \item We consider now the other inclusion.
    Let us first verify that
    $$M^{G_\delta}_{<\alpha}\subseteq M[G_\alpha].$$
    
    \item This is done by inductively showing
    $$M^{G_\delta}_\xi\subseteq M[G_\alpha]$$
    for all $\xi\in E\cap\alpha$.
    
    \item Assume that the statement is true all $\eta\in E\cap \xi$ and let us verify it for $\xi$.
    
    \item\label{635} By the assumption,
    $$M^{G_\delta}_{<\xi}=\bigcup_{\eta\in E\cap\xi}M^{G_\delta}_\eta\subseteq M[G_\alpha].$$
    
    \item If $\xi\not\in M^{G_\delta}_{<\xi}$, then
    $$M^{G_\delta}_{\xi}=M^{G_\delta}_{<\xi}\subseteq M[G_\alpha].$$
    Hence, let us assume $\xi\in M^{G_\delta}_{<\xi}$
    
    \item Let $x\in M^{G_\delta}_\xi$ be arbitrary.
    Then there exists $\dot{x}\in M^{G_\delta}_{<\xi}\cap V^{\P_\xi}$ such that $x=\dot{x}^{G_\xi}$.
    
    \item By \ref{635}, we have $\dot{x},\xi\in M[G_\alpha]$.
    
    \item Since also $G_\alpha\in M[G_\alpha]$, we conclude $G_\xi\in M[G_\alpha]$.
    
    \item Hence, $x=\dot{x}^{G_\xi}\in M[G_\alpha]$.
    
    \item This concludes the induction.
    Let us now go back to the main point, i.e. $M_\alpha [G_\delta]\subseteq M[G_\alpha]$.
    
    \item Let $x\in M_\alpha [G_\delta]$ be arbitrary.
    Then there exists $\dot{x}\in M^{G_\delta}_{<\alpha}\cap V^{\P_\alpha}$ such that $x=\dot{x}^{G_\alpha}$.
    
    \item By what we have just proved, we have $\dot{x}\in M [G_\alpha]$.
    
    \item Since $G_\alpha\in M[G_\alpha]$ as well, we conclude $x=\dot{x}^{G_\alpha}\in M[G_\alpha]$.
\end{proof0}

\item\notation Let $A$ be admissible, let $\Vec{\P}$ be an iteration with support $E$ that is point-wise definable over $A$, let $M\prec A$, let $\delta\in E$, let $G_{\delta}\leadsto V^{\P_{\delta}}$, and let $\alpha\in E\cap [0,\delta]$.
We shall henceforth write
$$M^{G_\alpha}:=M^{G_\delta}_\alpha,\quad M[G_\alpha]:=M_\alpha [G_\delta].$$
By the previous two propositions, this notation introduces no ambiguity.

\item\label{splitting iterations vm}\proposition Let $A$ be admissible, let $\Vec{\P}$ be an iteration with the support $E$ which is point-wise definable over $A$, let $M\prec A$, let $\alpha,\delta\in E$ with $\alpha\leq\delta$, and let $G_\delta\leadsto V^{\P_\delta}$.
Suppose that $\delta\in M^{G_\alpha}$.
Then 
$$M[G_\delta]=M^{G_\alpha}[G_\delta]\mbox{ and }M^{G_\delta}=(M^{G_\alpha})^{G_\delta}.$$
\begin{proof}
This is in fact obvious.
Namely, the second equality simply states that extending $M$ to $M^{G_\delta}$ according to the iteration $(\P_\xi : \xi\in E\cap [0,\delta])$ is the same as extending $M$ to $M^{G_\alpha}$ according to $(\P_\xi : \xi\in E\cap [0,\alpha])$ and then extending $M^{G_\alpha}$ to $(M^{G_\alpha})^{G_\delta}$ according to $(\P_\xi : \xi\in E\cap (\alpha,\delta])$.
An additional note is only required in the case $\alpha=\delta\in M^{G_\delta}$, where we apply \ref{338} of \ref{vm and genericity}.
For the first equality, we now have
$$M^{G_\alpha}[G_\delta]=(M^{G_\alpha})^{G_\delta}[G_\delta]=M^{G_\delta}[G_\delta]=(M^{G_\delta}_{<\delta})^{G_\delta}[G_\delta]=M^{G_\delta}_{<\delta}[G_\delta]=M[G_\delta],$$
where the first equality follows from \ref{vm and genericity}, the second follows from what we have just prove, the third by definition, the forth again by \ref{vm and genericity}, and the fifth by definition.
\end{proof}

\item\proposition Let $\gamma$ be such that $V_\gamma\models\mathsf{ZFC}^-$, let $A$ be admissible, let $M\prec A$ be a $\gamma$-strong virtual model, let $\Vec{\P}\in V_\gamma$ be an iteration with support $E$ which is point-wise definable over $A$, let $\alpha\in E$, and let $G_\alpha\leadsto V^{\P_\alpha}$.
Suppose that $M[G_\alpha]$ is defined.
Then $M[G_\alpha]$ is $\gamma$-strong in $V[G_\alpha]$.
\begin{proof}
This follows immediately from \ref{strength and forcing}.
\end{proof}

\item\label{reductions and iterated forcing 2}\proposition Let $\gamma$ be such that $V_\gamma\models\mathsf{ZFC}^-$, let $A$ be admissible, let $M,N\prec A$ be $\gamma$-strong virtual models, let $\Vec{\P}\in V_\gamma$ be an iteration with support $E$ which is point-wise definable over $A$, let $\alpha\in E$, and let $G_\alpha\leadsto V^{\P_\alpha}$.
Suppose that $M\cong_\gamma N$.
Then $M^{G_\alpha}\cong_\gamma N^{G_\alpha}$ and the isomorphism witnessing this fact extends the isomorphism witnessing $M\cong_\gamma N$.
\begin{proof}
This follows by induction on $\alpha$.
To handle the successor step, we simply apply \ref{cong-alpha and gen}.
Since the isomorphisms extend each other, the limit step is handled by simply taking the union of the previous isomorphisms.
\end{proof}

\item\label{cong-gamma and iterated ext}\proposition Let $\gamma$ be such that $V_\gamma\models\mathsf{ZFC}^-$, let $A$ be admissible, let $M,N\prec A$ be $\gamma$-strong virtual models, let $\Vec{\P}\in V_\gamma$ be an iteration with support $E$ which is point-wise definable over $A$, let $\alpha\in E$, and let $G_\alpha\leadsto V^{\P_\alpha}$.
Suppose that $M[G_\alpha]$ is defined.
Then $N[G_\alpha]$ is also defined and $M[G_\alpha]\cong_\gamma N[G_\alpha]$.
\begin{proof}
This follows from the previous proposition and \ref{cong-alpha and gen}.
\end{proof}

\end{enumerate}
\newpage

\section{Semi-proper Iteration}

\subsection{Setup for the Iteration}
\begin{enumerate}

\item\definition A \textit{blueprint for a semi-proper iteration} is a pair $(\V,\mathbf{Q})$ satisfying
\begin{parts}
    \item $\V=(V_\kappa,\in,U)$,
    
    \item $\kappa$ is inaccessible,
    
    \item $U\subseteq V_\kappa$,
    
    \item $\mathbf{Q}:\kappa\times V_\kappa\rightarrow V_\kappa$ is definable without parameters over $\V$,
    
    \item for all $\alpha<\kappa$ and for all $\P\in V_\kappa$, if $\P$ is a poset, then $\mathbf{Q}(\alpha,\P)\in V^\P$ and $\Vdash_\P$``$\mathbf{Q}(\alpha,\Q)$ is semi-proper".
\end{parts}

\item\declaration We fix a blueprint $(\V,\mathbf{Q})$ for a semi-proper iteration, where $\V=(V_\kappa,\in,U)$.
We will assume for the rest of the notes that the default language for virtual models is $\{\in,\dot{U}\}$, where $\dot{U}$ is a unary predicate.

\item\notation
\begin{parts}
    \item For $\alpha<\kappa$, we denote $\V_\alpha:=(V_\alpha,\in,U\cap V_\alpha).$
    
    \item We denote $\e:=\e_\V=\{\alpha<\kappa : \V_\alpha\prec\V\}$.
    
    \item $\e^+:=\{\alpha+1 : \alpha\in\e\}$
    
    \item $\e^*:=\e\cup\e^+$.
\end{parts}


\item\definition
\begin{parts}
    \item For $\alpha\in\e$, we say that a virtual model $M$ is \textit{$\alpha$-correct} if $\V_\alpha\prec\widehat{M}$.
    
    \item For $\alpha\in\e$, we define set
    $$\mathcal{C}_\alpha:=\{M : M\mbox{ is countable, }\alpha\mbox{-correct, and }\alpha\mbox{-generated virtual model}\}.$$
    
    \item For $S\subseteq\ord$, we define
    $$\mathcal{C}_S:=\bigcup_{\alpha\in S\cap\e}\mathcal{C}_\alpha.$$
\end{parts}

\item\remark Traditionally, forcing iterations are constructed as follows:
\begin{eqnarray*}
\P_0 & := & \{1\}\\
\P_{\alpha+1} & := & \P_\alpha*\dot{\Q}_\alpha\\
\P_\lambda & \subseteq & \{p : (\forall\alpha<\lambda)(p\upharpoonright\alpha\in\P_\alpha)\}
\end{eqnarray*}
($\lambda$ limit).
In the limit step, threads are required to satisfy some sort of a ``support condition".
In our iteration, we will double the steps: after adding a poset $\dot{\Q}_\alpha$, we will have an additional step of adding a ``scaffolding".
At the limit stages, the threading will now be controlled by that scaffolding.

\item\notation
\begin{parts}
    \item For convenience, the iteration will be indexed by $\e^*$ and the initial poset $\P_{\min\e}$ will not be equal to $\{1\}$.
    Nevertheless, poset $\P_{\min\e}$ will stil be trivial.
    
    \item Stages $\alpha\in\e^+$ correspond to adding a poset.
    
    \item Stages $\alpha$ where $\alpha$ is a successor point of $\e$ correspond to the operation of adding a scaffolding.
    
    \item Limit points of $\e$ are limit stages of the iteration.
    
    \item The canonical name for a $\P_\alpha$-generic will be denoted by $\G_\alpha$ (for $\alpha\in\e^*$).
    Accordingly, a $\P_\alpha$-generic will have denotation $G_\alpha$ and for $\beta\in\e^*\cap\alpha$, we will write $G_\beta:=G_\alpha\cap\P_\beta$.
\end{parts}

\item\remark The iteration $(\P_\alpha : \alpha\in\e^*)$ is defined by recursion.
The recursive step of the definition is described in the following subsection.
\end{enumerate}

\subsection{Recursive Step of the Definition}
\label{recursive step}
\begin{enumerate}
\item\declaration
\begin{parts}
    \item Let $\delta\in\e$ and let $\Vec{\P}=(\P_\alpha : \alpha\in\e^*\cap\delta)$ be a forcing iteration.
    
    \item Suppose that $\Vec{\P}$ is point-wise definable over $\V$.
    
    \item Suppose that for every $\alpha\in\e\cap\delta$ and every $p\in\P_\alpha$, it holds that that $p=(w_p,\M_p)$ for some finite partial function $w_p::\e\cap\alpha\rightarrow V_\kappa$ and some finite $\M_p\subseteq\mathcal{C}_{\leq\alpha}$.
    
    \item Let $(\dot{\Q}_\alpha : \alpha\in\e\cap\delta)$ be a sequence of names satisfying $\dot{\Q}_\alpha\in V^{\P_\alpha}$ and
    $$\P_\alpha\Vdash\dot{\Q}_\alpha\mbox{ is a semi-proper poset}$$
    for every $\alpha$.
\end{parts}


\item\definition Suppose that $\delta=\min\e$.
Then \textit{poset $\P_\delta$} consists of all pairs $p=(w_p,\M_p)$ where $w_p=\emptyset$ and $\M_p$ is a finite subset of $\mathcal{C}_{\delta}$.
\textit{The order $q\leq p$ in $\P_\delta$} holds if and only if $\M_q\supseteq\M_p$.

\item\definition Let $M\in\mathcal{C}_{\geq\delta}$, let $\alpha\in\e\cap\delta$, and let $G_\alpha\leadsto V^{\P_\alpha}$.
We say that \textit{$M^{G_\alpha}$ is active at $\delta$} in either of the following cases:
\begin{parts}
    \item there exists a predecessor $\gamma$ of ordinal $\delta$ inside $\e$ and $\gamma\in M^{G_\alpha}$;
    
    \item ordinal $\delta$ is a limit point of $\e$ and $\sup(M^{G_\alpha}\cap\e\cap\delta)=\delta$. 
\end{parts}

\item\definition Let $\M\in\mathcal{C}_{[0,\kappa]}$, let $\alpha\in\e\cap\delta$, and let $G_\alpha\leadsto V^{\P_\alpha}$.
We define
$$\M^\delta [G_\alpha]:=\{M\downarrow\delta : M\in\M,\, M^{G_\alpha}\mbox{ is active at }\delta\}.$$

\item\definition Let $\M\subseteq\mathcal{C}_\delta$, let $\alpha<\delta$, and let $G_\alpha\leadsto V^{\P_\alpha}$.
The set $\M$ is \textit{weak $\lhd_\delta$-chain w.r.t. $G_\alpha$} if for every $M,N\in\M$, it holds that:
\begin{parts}
    \item $\omega_1\cap M=\omega_1\cap N\implies M=N$;
    \item $\omega_1\cap M<\omega_1\cap N\implies M\lhd_\delta N^{G_\alpha}_\alpha$.
\end{parts}

\item\remark If in the second point the model $N^{G_\alpha}_\alpha$ is replaced by $N$, we get the usual notion of a $\lhd_\delta$-chain.

\item\remark Suppose that $\alpha\leq\beta<\delta$, that $G_\beta\leadsto V^{\P_\beta}$, and that $\M$ is a weak $\lhd_\delta$-chain w.r.t. $G_\alpha$.
Then $\M$ is a weak $\lhd_\delta$-chain w.r.t. $G_\beta$.

\item\label{retraction def}\notation Let $p=(w_p,\M_p)$ where $w_p$ is a partial function on $\e$ and $\M_p\subseteq\mathcal{C}_{[0,\kappa]}$ and let $\alpha\in\e$.
Then 
$$\M_p\downarrow\alpha:=\{M\downarrow\alpha : M\in\M_p\},$$
$$p\upharpoonright\alpha:=(w_p\upharpoonright\alpha, \M_p\downarrow\alpha),$$
$$p\upharpoonright(\alpha+1):=(w_p\upharpoonright(\alpha+1), \M_p\downarrow\alpha).$$

\item\definition Suppose that $\delta$ has a predecessor $\gamma$ inside $\e$.
The \textit{the poset $\P_\delta$} consists of all pairs $p=(w_p,\M_p)$ satisfying:
\begin{parts}[ref={[\arabic{section}.\arabic{subsection}.\arabic{enumi}\alph{partsi}]} ]
    \item object $w_p$ is a finite partial function $\e\cap\delta\rightarrow V_\kappa$;
    
    \item object $\M_p$ is a finite subset of $\mathcal{C}_{\leq\delta}$;
    
    \item $p\upharpoonright\alpha\in\P_\alpha$ for every $\alpha\in\e^*\cap\delta$;
    
    \item $p\upharpoonright(\gamma+1)\Vdash_{\P_{\gamma+1}}({\M}_p^\delta [\G_\gamma]\mbox{ is a weak }\lhd_\delta\mbox{-chain w.r.t. }\dot{G}_{\gamma+1})$
    
    \item\label{690} if $\gamma\in\dom(w_p)$, then
    $$p\upharpoonright\gamma\Vdash_{\P_\gamma}(\forall M\in{\M}^\delta_p[\G_\gamma])(w_p(\gamma)\mbox{ is semi-}M[\dot{G}_{\gamma}]^{\dot{\Q}_\gamma}\mbox{-generic}).$$
\end{parts}
\textit{The order $q\leq p$ in $\P_\delta$} holds if and only if:
\begin{parts}[resume]
    \item $q\upharpoonright\alpha\leq_{\P_\alpha} p\upharpoonright\alpha$ for every $\alpha\in\e^*\cap\delta$;
    
    \item $\M_q\cap\mathcal{C}_\delta\supseteq\M_p\cap\mathcal{C}_\delta$.
\end{parts}


\item\definition Suppose that $\delta$ is a limit point of $\e$.
Then \textit{the poset $\P_\delta$} consists of all pairs $p=(w_p,\M_p)$ satisfying:
\begin{parts}[ref={[\arabic{section}.\arabic{subsection}.\arabic{enumi}\alph{partsi}]} ]
    \item object $w_p$ is a finite partial function $\e\cap\delta\rightarrow V_\kappa$;
    
    \item object $\M_p$ is a finite subset of $\mathcal{C}_{\leq\delta}$;
    
    \item $p\upharpoonright\alpha\in\P_\alpha$ for every $\alpha\in\e^*\cap\delta$;
    
    \item there exists $\delta_0<\delta$ such that for every $\alpha\in\e\cap (\delta_0,\delta)$ we have that 
    $$p\upharpoonright(\alpha+1)\Vdash_{\P_{\alpha+1}}({\M}_p^\delta [\G_\alpha]\mbox{ is a weak }\lhd_\delta\mbox{-chain w.r.t. }\dot{G}_{{\alpha+1}}).$$
\end{parts}
\textit{The order $q\leq p$ in $\P_\delta$} holds if and only if:
\begin{parts}[resume]
    \item $q\upharpoonright\alpha\leq_{\P_\alpha} p\upharpoonright\alpha$ for every $\alpha\in\e^*\cap\delta$;
    
    \item $\M_q\cap\mathcal{C}_\delta\supseteq\M_p\cap\mathcal{C}_\delta$.
\end{parts}

\item\definition $\dot{\Q}_\delta:=\mathbf{Q}(\delta,\P_\delta)$

\item\definition \textit{Poset $\P_{\delta+1}$} consists of all pairs $(w_p,\M_p)$ satisfying:
\begin{parts}
    \item object $w_p$ is a finite partial function of the type $\e\cap (\delta+1)\rightarrow V_\kappa$;
    
    \item object $\M_p$ is a finite subset of $\mathcal{C}_{\leq\delta}$;
    
    \item $p\upharpoonright\delta\in\P_\delta$
    
    \item if $\delta\in\dom(w_p)$, then $w_p(\delta)$ is a canonical $\P_\delta$-name for an element of $\dot{\Q}_\delta$.
\end{parts}
\textit{The order $q\leq p$ in $\P_{\delta+1}$} holds if and only if:
\begin{parts}[resume]
    \item $q\upharpoonright\delta\leq_{\P_\delta} p\upharpoonright\delta$
    
    \item if $\delta\in\dom(w_p)$, then $\delta\in\dom(w_q)$ and $q\upharpoonright\delta\Vdash_{\P_\delta} w_q(\delta)\leq w_p(\delta).$
\end{parts}
\end{enumerate}

\subsection{Basic Properties}
\begin{enumerate}
\item\definition \textit{The semi-proper iteration given by the blueprint $(\V,\mathbf{Q})$} is obtained by recursively iterating the construction of Subsection \ref{recursive step}.

\item\definition Let $\P$ and $\Q$ be posets and let $\pi:\Q\rightarrow\P$.
Suppose that $\P$ is a suborder of $\Q$.
Then $\pi$ is said to be \textit{restriction of conditions} if
\begin{parts}
    \item $(\forall p\in\P)(\pi(p)=p)$,
    
    \item $(\forall p,q\in\Q)(p\leq q\implies \pi(p)\leq\pi(q))$,
    
    \item $(\forall p\in\P)(\forall q\in\Q)(p\leq \pi(q)\implies p\parallel q)$.
\end{parts}

\item\label{definition works}\proposition $(\P_\alpha : \alpha\in\e^*)$ is a correctly defined forcing iteration.
For elements $\alpha\leq\delta$ of $\e^*$, the mapping $\P_\delta\rightarrow\P_\alpha:p\mapsto p\upharpoonright\alpha$ is a restriction of conditions.
\begin{proof0}
    \item Let $\delta\in\e$ and let us assume recursively that $(\P_\alpha : \alpha\in\e^*\cap\delta)$ has been defined correctly and that it is an iteration.
    
    \item We can now define the poset $\P_\delta$.
    We need to verify $\P_\alpha$ is a complete subposet of $\P_\delta$ for every $\alpha\in\e^*\cap\delta$.
    
    \item It is immediate from the definitions that $\P_\alpha$ is a suborder of $\P_\delta$.
    It suffices for the conclusion to verify that the mapping $\P_\delta\rightarrow\P_\alpha : p\mapsto p\upharpoonright\alpha$ is a restriction of conditions, which is also obvious.
    
    \item We can now define the poset $\P_{\delta+1}$. 
    
    \item It is clear that $\P_\delta$ is a suborder of $\P_{\delta+1}$ and it is a matter of routine to verify that the mapping $\P_{\delta+1}\rightarrow\P_\delta:p\mapsto p\upharpoonright\delta$ is a restriction of conditions.
    Hence, $\P_\delta$ is a complete subposet of $\P_{\delta+1}$.
\end{proof0}

\item\definition $\P_\kappa:=\bigcup_{\alpha\in\e^*}\P_\alpha$

\item\proposition $(\P_\alpha : \alpha\in\e^*\cup\{\kappa\})$ is a forcing iteration.
Mapping $\P_\kappa\rightarrow\P_\alpha : p\mapsto p\upharpoonright\alpha$ is a restriction of conditions for all $\alpha\in\e^*$.\qed 

\item\label{easier definition}\proposition Let $\delta\in\e$.
\begin{itemize}
    \item Let $p=(w_p,\M_p)$.
    Then $p\in\P_\delta$ if and only if the following holds:
    \begin{parts}
        \item $w_p$ is a finite partial function from $\e\cap\delta$ into $V_\kappa$;
        
        \item $\M_p\in [\mathcal{C}_{\leq\delta}]^{<\omega}$
        
        \item for all $\gamma\in (\min\e,\delta]$, there exists $\beta\in\e\cap\gamma$ such that for all $\alpha\in\e\cap [\beta,\gamma)$,
        $$p\upharpoonright(\alpha+1)\Vdash_{\P_{\alpha+1}}(\M_p^\gamma [\G_\alpha]\mbox{ is a weak }\lhd_\gamma\mbox{-chain w.r.t. }\G_{\alpha+1});$$
        
        \item\label{974} for all $\beta\in\e\cap\delta$ and for $\gamma$ the successor of $\beta$ inside $\e$, 
        $$\beta\in\dom(w_p)\implies p\upharpoonright\beta\Vdash_{\P_\beta}(\forall M\in\M_p^\gamma [\G_\beta])(w_p(\beta)\mbox{ is semi-}M[\G_\beta]^{\dot{\Q}_\beta}\mbox{-generic}).$$
    \end{parts}
    
    \item Let $p,q\in\P_\delta$.
    Then $p\leq_{\P_\delta}q$ if and only if:
    \begin{parts}[resume]
        \item $\dom(w_p)\supseteq\dom(w_q)$
        
        \item for all $\alpha\in\dom(w_q)$,
        $$p\upharpoonright\alpha\Vdash_{\P_\alpha}w_p(\alpha)\leq_{\Q_\alpha}w_q(\alpha);$$
        
        \item for all $N\in\M_q$ and for $\gamma:=\sup\{\beta : V_\beta\subseteq \widehat{N}\}$, there exists $M\in\M_q$ such that $N=M\downarrow\gamma$.
    \end{parts}
\end{itemize}

\item\proposition Iteration $(\P_\alpha : \alpha\in\e^*)$ is point-wise definable over $\V$.
\qed 

\item\proposition Caridinality of $\P_\alpha$ is strictly less then $\beta$ for $\alpha\in\e^*$ and $\beta=\min(\e-(\alpha+1))$.
\qed 

\item\proposition $\P_{\delta+1}$ is forcing equivalent to $\P_\delta*\dot{\Q}_\delta$ for every $\delta\in\e$.
\qed 

\item\remark The previous equivalence is canonical.

\item\remark Suppose that $\beta<\gamma$ are elements of $\e^*$.
Recall that we have a canonical name $\P_\gamma/\P_\beta\in V^{\P_\beta}$ that satisfies
$$\Vdash_{\P_\beta} \P_\gamma/\P_\beta=\{q\in\P_\gamma : q\upharpoonright\beta\in\dot{G}_{\beta}\}.$$
    This is a name for a poset and the following is verified.
    \begin{parts}
        \item $\P_\beta*(\P_\gamma/\P_\beta)$ is forcing equivalent to $\P_\gamma$.
        
        \item For every dense subset $D$ of the poset $\P_\gamma$, the set $\{q\in D : q\upharpoonright\beta\in\dot{G}_{\P_\beta}\}$ is dense in the poset $\P_\gamma/\P_\beta$ inside the universe $V^{\P_\beta}$.
        
        \item If $G_\beta\leadsto V^{\P_\beta}$ and $H\leadsto V[G_\beta]^{(\P_\gamma/\P_\beta)^{G_\beta}}$, then $H$ is also a filter in $\P_\gamma$.
        If we want to think of it as such, we denote it by $G_\beta\cdot H$.
        Note that $G_\beta\cdot H$ corresponds to $G_\beta*H$ under the equivalence $\P_\gamma\simeq\P_\beta*(\P_\gamma/\P_\beta)$.
        In particular, $G_\beta\cdot H\leadsto V^{\P_\gamma}$ and $G_\beta\subseteq G_\beta\cdot H$.
        
    \end{parts}
    
    \item\label{meet}\definition Let $\alpha\leq\beta$ be elements of $\e^*$, let $p\in\P_\alpha$, and let $q\in\P_\beta$.
    Suppose that $p\leq q\upharpoonright\alpha$.
    Then we define $pq=(w,\M)$ as follows:
    $$w:=w_p\cup(w_q\upharpoonright(\alpha,\beta]),\quad \M:=\M_p\cup\M_q.$$
    
    \item\proposition Let $\alpha\leq\beta$ be elements of $\e^*$, let $p\in\P_\alpha$, and let $q\in\P_\beta$.
    Suppose that $p\leq q\upharpoonright\alpha$.
    Then $pq\in\P_\beta$, it satisfies $pq\leq p,q$ and $pq\upharpoonright\alpha=p$, and
    $$(\forall r\in\P_\beta)(r\leq p,q\iff r\leq pq).$$
    \qed 
    
\end{enumerate}

\subsection{Statement of Transfer Theorem}
\begin{enumerate}
\item\definition \Let
\begin{parts}
    \item $M$ : virtual model,
    
    \item $(\S_\alpha : \alpha\in E)$ : forcing iteration point-wise definable over $\widehat{M}$,
    
    \item $\gamma\in E$,
    
    \item $p\in\S_\gamma$.
\end{parts}
\Suppose there exists a beth-fixed point $\alpha$ such that $\Vec{\S}\in V_\alpha\subseteq\widehat{M}$.\\
\Then we say that \textit{$p$ is locally $M^{\S_\gamma}$-generic} if 
$$p\Vdash_{\S_\gamma}^V M^{\G_\gamma}=M^{\G_\gamma}_{<\gamma}.$$

\item\remark If $E=\{0\}$, then $p$ is locally $M^{\S_0}$-generic if and only if $p$ is $M^{\S_0}$-generic\footnote{in the sense of Remark 1.2-XI.f. of \cite{kasum2021iterating}}.
In general, ``$M^{\S_\gamma}$-generic" implies ``locally $M^{\S_\gamma}$-generic", but we will use this notion here to propagate semi-genericity through our iteration.

\item\textbf{Transfer Theorem.} \Let $\gamma\in\e$.
\Suppose for all $\alpha\in\e\cap\gamma$, for all $M\in\mathcal{C}_{>\alpha}$, for all $p\in\P_\alpha$ satisfying $M\downarrow\alpha\in\M_p$, $p\Vdash_{\P_\alpha}M^{\G_\alpha}\cap\omega_1^V=M\cap\omega_1^V$.
\Then for all $M\in\mathcal{C}_{>\gamma}$, every condition $p\in\P_\gamma$ satisfying $M\downarrow\gamma\in\M_p$ is locally $M^{\P_\gamma}$-generic.
\begin{proof}
The case $\gamma=\min\e$ is obvious since the poset $\P_\gamma$ is trivial.
The case where $\gamma$ is a successor point of $\e$ is proved in Subsubsection \ref{successor} and the case where $\gamma$ is a limit point of $\e$ is proved in Subsubsection \ref{limit}.
\end{proof}

\item\label{loc gen crit}\textbf{Local Genericity Criterion.}
\begin{itemize}[leftmargin=*, label={}]
    \item \Let
    \begin{itemize}[label=\textbullet]
        \item $(\S_\alpha : \alpha\in E)$ : forcing iteration,
        
        \item $M$ : virtual model,
        
        \item $\gamma\in E$,
        
        \item $p\in\S_\gamma$.
    \end{itemize}
    
    \item \Suppose 
    \begin{assum}
        \item there exists a beth-fixed point $\alpha$ such that $\Vec{\S}\in V_\alpha\subseteq\widehat{M}$,
        
        \item for all $q\leq_{\S_\gamma}p$, for all dense open subsets $D$ of $\S_\gamma$ satisfying $q\Vdash_{\S_\gamma}D\in M^{\G_\gamma}_{<\gamma}$, there exist $r\in D$ and $s\leq_{\S_\gamma}q,r$ such that $$s\Vdash_{\S_\gamma}\gamma\not\in M^{\G_\gamma}_{<\gamma}\vee r\in M^{\G_\gamma}_{<\gamma}.$$
    \end{assum}
    
    \item \Then $p$ is locally $M^{\S_\gamma}$-generic.
\end{itemize}
\begin{proof0}
    \item Let $p_0\leq p$ be arbitrary and let us show that there exists $p_1\leq p_0$ such that
    $$q\Vdash M^{\G_\gamma}=M^{\G_\gamma}_{<\gamma}.$$
    
    \item If $p_0\Vdash\gamma\not\in M^{\G_\gamma}_{<\gamma}$, then the conclusion follows by definition.
    Let us assume that $p_0\not\Vdash\gamma\not\in M^{\G_\gamma}_{<\gamma}$.
    
    \item Then there exists $p_1\leq p_0$ such that $p_1\Vdash\gamma\in M^{\G_\gamma}_{<\gamma}$.
    We claim that $p_1$ is as required.
    
    \item Let $\dot{\tau}\in V^{\S_\gamma}$ be such that $p_1\Vdash\dot{\tau}\in (M^{\G_\gamma}_{<\gamma})^{\S_\gamma}$.
    We want to show that $p_1\Vdash\dot{\tau}^{\G_\gamma}\in M^{\G_\gamma}_{<\gamma}$.
    
    \item Let $p_2\leq p_1$ be arbitrary.
    We want to find $s\leq p_2$ such that $s\Vdash\dot{\tau}^{\G_\gamma}\in M^{\G_\gamma}_{<\gamma}$.
    
    \item Note that $p_2\Vdash \dot{\tau}\in V^{\S_\gamma}$.
    This means that there exists $p_3\leq p_2$ and some $\sigma\in V^{\S_\gamma}$ such that $p_3\Vdash\dot{\tau}=\check{\sigma}$.
    
    \item Let $D$ be the set of all conditions in $\S_\gamma$ that decide the value of $\sigma$.
    Set $D$ is dense open in $\S_\gamma$.
    For $r\in D$, let $x_r$ be the value of $\sigma$ as decided by $r$.
    We have
    $$p_3\Vdash D, (x_r : r\in D)\in M^{G_\gamma}_{<\gamma}.$$
    
    \item There exist $r\in D$ and $s\leq p_3,r$ such that 
    $$s\Vdash r\in M^{G_\gamma}_{<\gamma}.$$
    
    \item Then
    $$s\Vdash\dot{\tau}^{\G_\gamma}=\check{\sigma}^{\G_\gamma}=x_s=x_r\in M^{G_\gamma}_{<\gamma},$$
    the last fact being due to $s\Vdash (x_r)_{r\in D},r\in M^{G_\gamma}_{<\gamma}$.
\end{proof0}
\item\remark In the case $E=\{\emptyset\}$, the above criterion reduces to the usual genericity criterion (see for example Lemma 1.2-IV of \cite{kasum2022proper}).
\end{enumerate}

\subsection{Some Lemmas for Transfer Theorem}
\begin{enumerate}
\item\label{adding a side condition}\textbf{Lemma.} Let $\beta\in\e$, let $N\in\mathcal{C}_{\geq\beta}$, and let $p\in\P_\beta\cap N$.
Then there exists $q\leq_{\P_\beta}p$ such that $\dom(w_q)=\dom(w_p)$ and $\M_q=\M_p\cup\{N\downarrow\beta\}$.
\begin{proof0}
    \item Let $\M_q:=\M_p\cup\{N\downarrow\beta\}$.
    We need to define $w_q$ in such a way as to ensure $q:=(w_q,\M_q)\in\P_\beta$ and $q\leq p$.
    
    \item Note that for all $P\in\M_p$, we have $P\cap\omega_1^V<N\cap\omega_1^V$ and $P\lhd_\beta N$.
    Hence, we will be done with the argument if we construct $w_q$ satisfying:
    \begin{parts}
        \item $\dom(w_q)=\dom(w_p)=:d$
        
        \item $q\upharpoonright\alpha\Vdash_{\P_\alpha}w_q(\alpha)\leq w_p(\alpha)$ for all $\alpha\in d$;
        
        \item $q\upharpoonright\alpha\Vdash_{\P_\alpha}\left(\alpha\in N^{\G_\alpha}\implies (w_p(\alpha)\mbox{ is semi-}N[\G_\alpha]^{\dot{\Q}_\alpha}\mbox{-generic})\right)$ for all $\alpha\in d$.
    \end{parts}
    (Cf. Proposition \ref{easier definition}.)
    
    \item $w_q(\alpha)$ is defined by recursion on $\alpha\in d$.
    Suppose that $w_q\upharpoonright\alpha$ has been defined and let us define $w_q(\alpha)$.
    
    \item Let $G_\alpha\leadsto V^{\P_\alpha}$ be an arbitrary generic containing $q\upharpoonright\alpha$.
    
    \item Suppose that $\alpha\in N^{G_\alpha}$.
    Then $N[G_\alpha]$ is defined.
    
    \item Since $\Q_\alpha:=\dot{\Q}_\alpha^{G_\alpha}\in N[G_\alpha]$ is semi-proper and $w_p(\alpha)^{G_\alpha}\in N[G_\alpha]$, there exists $w\leq_{\Q_\alpha}w_p(\alpha)^{G_\alpha}$ which is semi-$N[G_\alpha]^{\Q_\alpha}$-generic.
    
    \item Since $G_\alpha$ was arbitrary containing $q\upharpoonright\alpha$, we can find a name $\dot{w}\in V^{\P_\alpha}$ such that
    $$q\upharpoonright\alpha\Vdash_{\P_\alpha}\dot{w}\leq_{\dot{\Q}_\alpha}w_p(\alpha)\wedge\left(\alpha\in N^{\G_\alpha}\implies(\dot{w}\mbox{ is semi-}N^{G_\alpha}\mbox{-generic})\right).$$
    
    \item We can now set $w_q(\alpha):=\dot{w}$.
\end{proof0}

\item\label{1325}\textbf{Lemma.} \underline{Let}
\begin{itemize}
    \item $\mu,\nu\in\e$ : $\nu$ is the successor of $\mu$ in $\e$,
    
    \item $p\in\P_\nu$,
    
    \item $G_\mu\leadsto V^{\P_\mu}$ : $p\upharpoonright\mu\in G_\mu$,
    
    \item $M\in\M_p^\nu [G_\mu]$,
    
    \item $\Q_\mu:=\dot{\Q}_\mu^{G_\mu}$,
    
    \item $u\in\Q_\mu\cap M[G_\mu]$.
\end{itemize}
\underline{Then} there exists $v\in\Q_\mu$ such that $v\leq u$ and that for all $N\in\M_p^\nu [G_\mu]$ satisfying $N\cap\omega_1^V\geq M\cap\omega_1^V$, it holds that $v$ is semi-$N[G_\mu]^{\Q_\mu}$-generic.
\begin{proof0}
    \item Let
    \begin{itemize}
        \item $\{M_i : i<n\}:=\M_p^\nu [G_\mu]$ : for all $i<j<n$, $M_i\cap\omega_1^V<M_j\cap\omega_1^V$,
        
        \item $k<n$ : $M_k=M$,
        
        \item $\lambda:=|\Q_\mu|$.
    \end{itemize}
    
    \item For all $i<n$, we have $\Q_\mu\in M_i[G_\mu]$ and $\Q_\mu$ is semi-proper (this follows from activity).
    
    \item We may assume w.l.o.g. that $\Q_\mu=\lambda$.
    
    \item Consider $i<n$.
    We let $N_i:=M_i[G_\mu]\cap H((2^\lambda)^+)\prec H((2^\lambda)^+)$.
    Note that $w\in\Q_\mu$ is semi-$M_i[G_\mu]^{\Q_\mu}$-generic if and only if it is semi-$N_i^{\Q_\mu}$-generic.
    
    \item Hence, it suffices to find $v\in\Q_\mu$ such that $v\leq u$ and which is semi-$N_i^{\Q_\mu}$-generic for all $i\in [k,n)$.
    
    \item Since $M_i\lhd_\nu M_j[G_\mu]$ for all $i<j<n$, Proposition \ref{alpha-in to in for H} implies that $N_i\in N_j$ for all $i<j<n$.
    
    \item By recursion on $i\in [k,n)$, we construct a sequence $(u_i : k\leq i<n)$ of $\Q_\mu$-conditions satisfying $u_k:=u\in N_k$ and satisfying for all $i\in (k,n)$ that $u_i\in N_i$, $u_i\leq_{\Q_\mu}u_{i-1}$, and $u_i$ is semi-$N_{i-1}^{\Q_\mu}$-generic.
    
    \item There exists $v\leq u_{n-1}$ such that $v$ is semi-$N_{n-1}^{\Q_\mu}$-generic.
\end{proof0}

\item\label{1372}\textbf{Lemma.} \underline{Let}
\begin{itemize}
    \item $\mu,\nu\in\e$ : $\nu$ is the successor of $\mu$ in $\e$,
    
    \item $p\in\P_\nu$,
    
    \item $M\in\M_p$ : $p\upharpoonright\mu\Vdash_{\P_\mu}(M^{\G_\mu}$ is active at $\nu)$,
    
    \item $\dot{u}\in V^{\P_\mu}$ : $p\upharpoonright\mu\Vdash_{\P_\mu}\dot{u}\in\dot{\Q}_\mu\cap M[\G_\mu]$.
\end{itemize}
\underline{Then} there exists a canonical $\P_\mu$-name $\dot{v}$ for an element of $\dot{\Q}_\mu$ such that $p\upharpoonright\mu$ forces that ``$\dot{v}\leq\dot{u}$ and that for all $N\in\M_p^\nu [\G_\mu]$ satisfying $N\cap\omega_1^V\geq M\cap\omega_1^V$, it holds that $\dot{v}$ is semi-$N[\G_\mu]^{\dot{\Q}_\mu}$-generic.
\begin{proof}
This is immediate from the previous lemma.
\end{proof}

\item\label{adding to the working part}\lemma Let $\gamma\in\e$, let $p\in\P_\gamma$, and let $\beta\in\e\cap\gamma$.
Then there exists $p_\beta\leq_{\P_\gamma} p$ such that $\beta\in\dom(w_{p_\beta})$.
\begin{proof0}
    \item We may assume w.l.o.g. that $(\beta,\gamma)\cap\e=\emptyset$.
    
    \item Let $G_\beta\leadsto V^{\P_\beta}$ be arbitrary containing $p\upharpoonright\beta$, let $M\in\M_p^\gamma [G_\beta]$ be such that $M\cap\omega_1^V$ is minimal, and let $\Q_\beta:=\dot{\Q}_\beta^{G_\beta}$.

    \item We may apply Lemma \ref{1325} to
    \begin{parts}
        \item $\underaccent{\tilde}{\mu}:=\beta$, $\underaccent{\tilde}{\nu}:=\gamma$,
        
        \item $\underaccent{\tilde}{p}:=p\upharpoonright\beta$,
        
        \item $\underaccent{\tilde}{G}_\mu:=G_\beta$,
        
        \item $\underaccent{\tilde}{M}:=M$,
        
        \item $\underaccent{\tilde}{u}:=1_{\Q_\beta}$,
    \end{parts}
    which yields
    $$(\exists v\in\Q_\beta)(\forall N\in\M_p^\gamma [G_\beta])(v\mbox{ is semi-}N[G_\beta]^{\Q_\beta}\mbox{-generic}).$$
    
    \item Since $G_\beta$ was arbitrary, we can find a canonical $\P_\beta$ name $\dot{v}$ for an element of $\dot{\Q}_\beta$ such that 
    $$p\upharpoonright\beta\Vdash_{\P_\beta}^V(\forall N\in\M_p^\gamma [\G_\beta])(\dot{v}\mbox{ is semi-}N[\G_\beta]^{\dot{\Q}_\beta}\mbox{-generic}).$$
    
    \item Let $w_{p_\beta}:=w_p\cup\{(\beta,\dot{v})\}$, $\M_{p_\beta}:=\M_p$, and $p_\beta:=(w_{p_\beta},\M_{p_\beta})$.
    We see that $p_\beta$ is as required.
\end{proof0}

\item\label{collapsing ordinals}\textbf{Lemma.} Let $\beta\in\e-\{\min\e\}$ and let $\xi<\beta$.
Then $\Vdash_{\P_\beta}|\xi|\leq\omega_1$.
\begin{proof0}
    \item It suffices to consider the case when $\beta$ has a predecessor $\alpha$ inside $\e$.
    
    \item Let $\dot{\M}$ be a $\P_\beta$-name satisfying
    $$\Vdash_{\P_\beta}\dot{\M}=\bigcup_{p\in\G_\beta}\M_p^\beta [\G_\alpha].$$
    
    \item We want to show that:
    \begin{parts}[ref=\arabic{pfenumi}$^\circ$\alph{partsi}]
        \item\label{1092} $\Vdash_{\P_\beta}\xi\subseteq\bigcup\dot{\M}$
        
        \item\label{1100} $\Vdash_{\P_\beta}|\dot{\M}|\leq\omega_1$.
    \end{parts}
    This suffices for the conclusion.
    
    \item Let us first verify \ref{1092}.
    Let $\eta<\xi$ and let $D:=\{p\in\P_\beta : (\exists M\in\M_p)(\eta,\alpha\in M)\}$.
    We are done if we show that $D$ is dense in $\P_\beta$.
    
    \item Let $p_0\in\P_\beta$ be arbitrary.
    Then there exists $N\prec\V$ such that $\alpha,\beta,\eta,p_0\in N$.
    
    \item By Lemma \ref{adding a side condition}, there exists $p\leq_{\P_\beta}p_0$ such that $N\downarrow\beta\in\M_p$.
    
    \item Since $\eta,\alpha\in N\downarrow\beta$, we conclude $p\in D$.
    
    \item Let us now verify \ref{1100}.
    Let $G_\beta\leadsto V^{\P_\beta}$ be arbitrary and let us work inside $V[G_\beta]$.
    
    \item Let $\M:=\dot{\M}^{G_\beta}\subseteq\mathcal{C}_\beta$.
    It suffices to show that the mapping
    $$\M\rightarrow\omega_1^V:M\mapsto M\cap\omega_1^V$$
    is an injection.
    
    \item Let $M,N\in\M$ be arbitrary satisfying $M\cap\omega_1^V=N\cap\omega_1^V$.
    We want to show that $M=N$.
    
    \item We have that $\alpha\in M\cap N$ and there exist $p,q\in G_\beta$ such that $M\in\M_p^\beta [G_\alpha]$ and $N\in\M_q^\beta [G_\alpha]$.
    
    \item Since $p,q\in G_\beta$, there exists $r\in G_\beta$ such that $r\leq_{\P_\beta}p,q$.
    We have that $M,N\in\M_r^\beta [G_\alpha]$.
    
    \item Since $r\in G_\beta$, we have that $\M_r^\beta [G_\alpha]$ is a weak $\lhd_\alpha$-chain w.r.t. $G_{\alpha+1}$.
    In particular, the fact that $M\cap\omega_1^V=N\cap\omega_1^V$ implies that $M=N$.
\end{proof0}

\item\label{1027}\textbf{Lemma.} Let $M$ be a virtual model, let $\delta$ be an inaccessible of $M$, let $\P\in M\cap V_\delta$ be a poset, and let $G\leadsto \widehat{M}^\P$.
Then
$$\sup(M[G]\cap\delta)=\sup(M\cap\delta).$$
\begin{proof0}
    
    
    

    \item Let $\tau\in M^\P$ be such that $\tau_G<\delta$.
    We want to show that there exists some $\beta\in M\cap\delta$ such that $\tau_G<\beta$.

    \item We may assume w.l.o.g. that $\Vdash_\P^V\tau<\delta$.

    \item There exists a maximal antichain $A\in M$ of $\P$ such that for all $p\in A$, there exists $\alpha_p<\delta$ such that $p\Vdash_\P^V\tau=\alpha_p$.

    \item Note that $(\alpha_p : p\in A)\in M$.
    Since $\delta$ is inaccessible in $M$, we conclude that $\beta:=\sup_{p\in A}\alpha_p\in M\cap\delta$.

    \item It is now clear that $\beta$ is as required.
\end{proof0}

\item\label{bounding activity}\lemma 
\begin{itemize}[leftmargin=*]
    \item[] \Let 
    \begin{itemize}[label=\textbullet]
        \item $\gamma$ : limit point of $\e$,
        
        \item $M\in\mathcal{C}_{[0,\kappa]}$,
        
        \item $p\in\P_\gamma$ : $M\downarrow\gamma\in\M_p$.
    \end{itemize}
    
    \item[] \Suppose for all $\alpha\in\e\cap\gamma$, condition $p\upharpoonright\alpha\Vdash_{\P_\alpha}M^{\G_\alpha}\cap\omega_1^V=M\cap\omega_1^V$.
    
    \item[] \Then there exist ordinals $\beta<\gamma^*\leq\gamma$ and a condition $q\leq_{\P_\gamma} p$ such that $\M_q\cap\mathcal{C}_\gamma=\M_p\cap\mathcal{C}_\gamma$ and
    $$(\forall \alpha\in\e\cap (\beta,\gamma))(q\upharpoonright\alpha\Vdash_{\P_\alpha}\sup(M^{\G_\alpha}\cap\gamma)=\gamma^*).$$
\end{itemize}
\begin{proof0}
    \item If $M\cap\ord\subseteq\gamma$, then we may set $\gamma^*:=\sup(M\cap\ord)$, $\beta:=0$, and $q:=p$.
    Hence, let us assume that there exists $\xi\in M$ with $\xi\geq\gamma$.
    
    \item Let $\eta$ be the least ordinal $\eta\geq\gamma$ satisfying
    $$(\exists\beta\in\e^*\cap\gamma)(\exists r\leq_{\P_\beta}p\upharpoonright\beta)(r\Vdash_{\P_\beta}\eta\in M^{\G_\beta}).$$
    
    \item Let $\beta_0\in\e^*\cap\gamma$ be the least ordinal satisfying
    $$(\exists r\leq_{\P_{\beta_0}}p\upharpoonright\beta_0)(r\Vdash_{\P_{\beta_0}}\eta\in M^{\G_{\beta_0}}).$$
    Let $r$ be a witness for the last formula.
    
    \item\label{1407} By lemma \ref{722}, we have that $\eta\in\e$.
    
    \item\claim $r\Vdash_{\P_{\beta_0}}\beta_0\in M^{\G_{\beta_0}}_{<\beta_0}$
    \begin{proof0}
        \item Assume otherwise.
        Then there exists $r'\leq_{\P_{\beta_0}}r$ such that $r'\Vdash\beta_0\not\in M^{\G_{\beta_0}}_{<\beta_0}$.
        
        \item This means that $r'\Vdash \eta\in M^{\G_{\beta_0}}=M^{\G_{\beta_0}}_{<\beta_0}$.
        
        \item Hence, there are $r''\leq_{\P_{\beta_0}}r'$ and $\beta_1\in\e^*\cap\beta_0$ such that $r''\Vdash\eta\in M^{\G_{\beta_1}}$.
        This contradicts the minimality of $\beta_0$.
    \end{proof0}
    
    \item Hence, $r\Vdash_{\P_{\beta_0}}(M[\G_{\beta_0}]$ is defined$)$.
    
    \item Let $s\leq_{\P_{\beta_0}}r$ and let $\lambda\in\ord$ be such that $s\Vdash_{\P_{\beta_0}}\cof^{M[\G_{\beta_0}]}(\eta)=\lambda$.
    
    \item\claim There exists $\beta\in\e\cap(\beta_0,\gamma)$ and $t\leq_{\P_\beta}s$ such that
    $$t\Vdash_{\P_\beta}\beta\in M^{\G_\beta}_{<\beta}\wedge\cof^{M[G_\beta]}(\eta)\in\{\omega,\omega_1,\eta\}.$$
    \begin{proof0}
        \item If $\lambda=\eta$, it suffices to take $\beta:=\min(\e-\beta_0)$ and $t:=s$.
        Hence, let us assume that $\lambda<\eta$.
        
        \item Note that $s\Vdash_{\P_{\beta_0}}\lambda\in M^{\G_{\beta_0}}$.
        By the choice of $\eta$, we then have that $\lambda<\gamma$.
        
        \item Let $G_{\beta_0}\leadsto V^{\P_{\beta_0}}$ be an arbitrary generic containing $s$ and let us work in $V[G_{\beta_0}]$.
        
        \item Since $\eta\in\e$, the fact that $\e\cap(\lambda,\eta)\not=\emptyset$ can be expressed as saying that there exists $\beta\in(\lambda,\eta)$ such that $\widehat{M}\upharpoonright\beta\prec\widehat{M}\upharpoonright\eta$.
        Since this is a first order fact of $\widehat{M}$ with parameters from $M^{G_{\beta_0}}\prec\widehat{M}$, there exists such an ordinal $\beta$ in $M^{G_{\beta_0}}$.
        
        \item By the choice of $\eta$, we have $\beta\in\e\cap(\lambda,\gamma)$.
        
        \item Let $G_\beta\leadsto V^{\P_\beta}$ extend $G_{\beta_0}$ and let us work inside $V[G_\beta]$.
        
        \item We have $\beta\in M^{G_{\beta_0}}\subseteq M^{G_\beta}_{<\beta}$.
        
        \item By Lemma \ref{collapsing ordinals}, we have $|\lambda|\leq\omega_1$ and consequently $\cof^{M[G_\beta]}(\eta)\leq\omega_1$.
        
        \item Condition $t$ is now obtained by Forcing Theorem.
    \end{proof0}
    
    \item Up to strengthening $t$, we may assume that there exists $\gamma^*$ such that $t\Vdash_{\P_\beta}\gamma^*=\sup(M[\G_{\beta}]\cap\gamma)$.
        
    Let $q:=tp\in\P_\gamma$ and let us show that $\beta$, $\gamma^*$, and $q$ are as required.
    It is immediate that $\M_q\cap\mathcal{C}_\gamma=\M_p\cap\mathcal{C}_\gamma$.
    
    \item Since $t\Vdash_{\P_{\beta}}\beta\in M[\G_\beta]\cap\gamma\wedge \gamma^*=\sup(M[\G_\beta]\cap\gamma)$, we conclude $\beta<\gamma^*\leq\gamma$.
    
    
        
        
        
        
        
        
    
    \item\claim For all $\alpha\in\e\cap (\beta,\gamma)$, we have that $q\upharpoonright\alpha\Vdash_{\P_\alpha}\sup(M^{\G_\alpha}\cap\gamma)=\gamma^*.$
    \begin{proof0}
        \item Let $G_\alpha\leadsto V^{\P_\alpha}$ be an arbitrary generic containing $q\upharpoonright\alpha$ and let us work inside $V[G_\alpha]$.
        We will distinguish between two cases, depending on whether $\cof^{M[G_\beta]}(\eta)\leq\omega_1^{V[G_\beta]}$ or $\cof^{M[G_\beta]}(\eta)=\eta$.
        
        
        
        
        
        \item \textit{Suppose first that $\cof^{M[G_\beta]}(\eta)\leq\omega_1^{V[G_\beta]}$.}
        Then there exists $f\in M[G_\beta]$
        such that 
        $$f:\omega_1^{V[G_\beta]}\rightarrow\eta$$
        and 
        $$M[G_\beta]\models(f\mbox{ is cofinal in }\eta).$$
        
        \item By the assumption of the lemma, we have that
        $$M[G_\beta]\cap\omega_1^V=M\cap\omega_1^V.$$
        In particular, $\omega_1^V\not\subseteq M[G_\beta]$, which means that $\omega_1^V$ is not countable in $V[G_\beta]$.
        
        \item Hence, $\omega_1^V=\omega_1^{V[G_\beta]}$ and $f:\omega_1^V\rightarrow\eta$ is such that 
        $$\sup(f[M\cap\omega_1^V])=\sup(M^{G_\beta}\cap\eta).$$
        
        \item\label{1143} This further implies that
        $$\sup(f[M\cap\omega_1^V])=\sup(M^{G_\beta}\cap\gamma)=\gamma^*.$$
        
        \item Observe that $M[G_\beta]\prec\widehat{M}[G_\beta]$, that $\widehat{M}[G_\alpha]$ is a generic extension of $\widehat{M}[G_\beta]$, and that $M[G_\alpha]\prec\widehat{M}[G_\alpha]$.
        This implies that 
        $$M[G_\alpha]\models(f:\omega_1^V\rightarrow\eta\mbox{ cofinally}).$$
        
        \item We can now compute as follows:
        $$\sup(M^{G_\alpha}\cap\gamma)=\sup(M^{G_\alpha}\cap\eta)=\sup(f[M^{G_\alpha}\cap\omega_1^V])=\sup(f[M\cap\omega_1^V])=\gamma^*,$$
        where the first equality follows by the choice of $\eta$, the second one from the previous point, 
        the third one from the hypothesis that
        $$M^{G_\alpha}\cap\omega_1^V=M\cap\omega_1^V,$$
        and the fourth one from \ref{1143}.
        
        \item \textit{Suppose now that $\cof^{M[G_\beta]}(\eta)=\eta$.}
        
        \item Since $M[G_\beta]\cap\eta=M^{G_\beta}\cap\eta=M^{G_\beta}\cap\gamma$, we see that $\eta$ is a strong limit in $M[G_\beta]$.
        Hence,
        $$M[G_\beta]\models(\eta\mbox{ is inaccessible}).$$
        
        \item Since posets $(\P_\xi : \xi\in\e^*\cap\gamma)$ are of rank $<\gamma\leq\eta$, the iterative application of Lemma \ref{1027} yields that
        $$\sup(M^{G_\xi}\cap\eta)=\sup(M^{G_\beta}\cap\eta),$$
        for all $\xi\in\e^*\cap [\beta,\alpha]$.
        
        \item In particular, we have that $\sup(M^{G_\alpha}\cap\eta)=\gamma^*$, which further implies
        $$\sup(M^{G_\alpha}\cap\gamma)=\gamma^*.$$
    \end{proof0}
    
    
    \item This completes the proof of the lemma.
\end{proof0}

\item\label{lifting vm}\lemma \Let
\begin{parts}
    \item $M$ : virtual model,
    
    \item $\gamma\in\e_M$ : limit point of $\e_{\widehat{M}}$,
    
    \item $\gamma^*:=\sup(M\cap\gamma)$,
    
    \item $N$ : virtual model satisfying $N\lhd_{\gamma^*}M$.
\end{parts}
\Then there exists a unique $\gamma$-generated model $N^+\in M$ such that $N\cong_{\gamma^*} N^+$.
\underline{Moreover}, if $\widehat{M}\upharpoonright\gamma^*\prec\widehat{N}$, then $\widehat{M}\upharpoonright\gamma\prec\widehat{N^+}$.
\begin{proof0}
    \item To establish \underline{existence} of such a $N^+$, note that by definition there exists $P\in M$ such that $N\cong_{\gamma^*}P$.
    Then $N^+:=P\downarrow\gamma\in M$ is as required.
    
    \item Let us verify \underline{uniqueness}.
    Let $N'\in M$ be a $\gamma$-generated virtual model such that $N\cong_{\gamma^*}N'$.
    It suffices to show $N'\cong_{\gamma}N^+$.
    
    \item Let $\alpha\in\e_M\cap\gamma$ be arbitrary.
    Then $\alpha<\gamma^*$ and consequently $N'\cong_\alpha N\cong_\alpha N^+$.
    
    \item Hence,
    $$M\models(\forall\alpha<\gamma)(V_\alpha\prec (V_\gamma,\in,U\cap V_\gamma)\implies N'\cong_\alpha N^+).$$
    
    \item By elementarity, the same statement is true in $\widehat{M}$.
    This means that for all $\alpha\in\e_{\widehat{M}}\cap\gamma$, we have that $N'\cong_\alpha N^+$.
    
    \item For $\alpha\in\e_{\widehat{M}}\cap\gamma$, let $f_\alpha:\hull(N',V_\alpha)\cong\hull(N^+,V_\alpha)$ be the unique isomorphism witnessing $N'\cong_\alpha N^+$.
    For $\alpha_0<\alpha_1$, we have $f_{\alpha_0}\subseteq f_{\alpha_1}$.
    
    \item It is now easily seen that 
    $$f_\gamma:=\bigcup_{\alpha\in\e_{\widehat{M}}\cap\gamma}f_\alpha:\hull(N',V_\gamma)\rightarrow\hull(N^+,V_\gamma)$$
    witnesses $\hull(N',V_\gamma)\cong_\gamma\hull(N^+,V_\gamma)$.
    
    \item Let us establish \underline{moreover part}.
    Let $\alpha\in\e_{M}\cap\gamma$ be arbitrary.
    We have that $\alpha<\gamma^*$.
    
    \item This implies that $\widehat{M}\upharpoonright\alpha\prec\widehat{N}$
    and consequently
    $\widehat{M}\upharpoonright\alpha=\hull(N,V_{\gamma^*})\cap V_\alpha\prec\hull(N,V_{\gamma^*})$.
    
    \item Since $N\cong_{\gamma^*}N^+$, we conclude that
    $$\widehat{M}\upharpoonright\alpha=\hull(N^+,V_{\gamma^*})\cap V_\alpha\prec\hull(N^+,V_{\gamma^*})\prec\widehat{N^+}$$
    and consequently $M\models V_\alpha\prec \widehat{N^+}$.
    
    \item Hence,
    $$M\models(\forall\alpha<\gamma)(V_\alpha\prec (V_\gamma,\in,U\cap V_\gamma)\implies V_\alpha\prec\widehat{N^+}).$$
    
    \item Since $M\prec\widehat{M}$, we conclude that for all $\alpha\in\e_{\widehat{M}}\cap\gamma$, we have $\widehat{M}\upharpoonright\alpha\prec\widehat{N^+}$.
    This suffices for the conclusion.
\end{proof0}
\end{enumerate}

\subsection{Proof of Transfer Theorem}

\subsubsection{Successor Case}
\label{successor}
\begin{pfenum}
\item Let $\gamma\in\e$ be the successor of some $\beta$ inside $\e$.
We are assuming that
\begin{equation}
\label{1546}
    (\forall N\in\mathcal{C}_{>\beta})(\forall q\in\P_\beta)(N\downarrow\beta\in\M_q\implies(q\mbox{ is semi-}N^{\P_\beta}\mbox{-generic})).
\end{equation}

Let $M\in\mathcal{C}_{>\gamma}$ and let $p\in\P_\gamma$ satisfying $M\downarrow\gamma\in\M_p$.
We want to show that $p$ is locally $M^{\P_\gamma}$-generic.

\item\label{1645} We will use Local Genericity Criterion \ref{loc gen crit}.
Let $D$ be a dense open subset of $\P_\gamma$ and let $q\leq_{\P_\gamma}p$ be arbitrary satisfying
$$q\Vdash_{\P_\gamma}D\in M^{\G_\gamma}_{<\gamma}.$$
We want to find $r\in D$ and $s\leq_{\P_\gamma}q,r$ such that 
$$s\Vdash_{\P_\gamma}\gamma\not\in M^{\G_\gamma}_{<\gamma}\vee r\in M^{\G_\gamma}_{<\gamma}.$$

\item Up to strengthening $q$, we may assume that $q\in D$.

\item If $q\Vdash_{\P_\gamma}\gamma\not\in M^{\G_\gamma}_{<\gamma}$, we can take $s:=r:=q$.
Hence, let us assume that $q\not\Vdash_{\P_\gamma}\gamma\not\in M^{\G_\gamma}_{<\gamma}$.

\item Up to strengthening $q$, we may assume that $q\Vdash_{\P_\gamma}\gamma\in M^{\G_\gamma}_{<\gamma}$.
Since $\beta$ is definable from $\gamma$ over $\widehat{M}$, we have that $q\Vdash_{\P_\gamma}\beta\in M^{\G_\gamma}_{<\gamma}$.
This implies that $q\Vdash_{\P_\beta}\beta\in M^{\G_\beta}$ and $q\Vdash_{\P_\beta}M\downarrow\gamma\in\M_q^\gamma [\G_\beta]$.

\item By Lemma \ref{adding to the working part}, we may strengthen $q$ to ensure $\beta\in\dom(w_q)$.

\item\label{1571}\claim $q\upharpoonright(\beta+1)\Vdash_{\P_{\beta+1}}M^{\G_{\beta+1}}\cap\omega_1^V=M\cap\omega_1^V$
\begin{proof0}
    \item By (\ref{1546}), we have that
    $$q\upharpoonright\beta\Vdash_{\P_\beta}M^{\G_{\beta}}\cap\omega_1^V=M\cap\omega_1^V.$$
    
    \item Let $G_{\beta+1}$ be a generic containing $q\upharpoonright(\beta+1)$ and let us verify that $M^{G_{\beta+1}}\cap\omega_1^V=M\cap\omega_1^V$.
    Let $\Q_\beta:=\dot{\Q}_\beta^{G_\beta}$.
    
    \item Since $\beta\in\dom(w_q)$ and $M\downarrow\gamma\in\M_q^\gamma [G_\beta]$, we have that $w_q(\beta)^{G_\beta}$ is semi-$(M\downarrow\gamma)[G_\beta]^{\Q_\beta}$-generic in $V[G_\beta]$.
    
    \item Hence,
    $$(M\downarrow\gamma)[G_{\beta+1}]\cap\omega_1^{V[G_\beta]}=(M\downarrow\gamma)[G_\beta]\cap\omega_1^{V[G_\beta]}.$$
    
    \item This shows that
    $$M^{G_{\beta+1}}\cap\omega_1^V=M^{G_\beta}\cap\omega_1^V=M\cap\omega_1^V.$$
\end{proof0}

\item\label{1566} We may assume that there exists $\M\in [\mathcal{C}_\gamma]^{<\omega}$ such that 
$${q}\upharpoonright\beta\Vdash_{\P_\beta}\M_q^\gamma [\G_\beta]=\M.$$
Namely, there exists $\bar{q}\leq_{\P_\beta}q\upharpoonright\beta$ such that $\bar{q}$ decides $\M_q^\gamma [\G_\beta]$.
Since $\M_{\bar{q}q}\cap\mathcal{C}_\gamma=\M_q\cap\mathcal{C}_\gamma$, we may replace $q$ by $\bar{q}q$.

\item Note that $M\downarrow\gamma\in\M$ and
$${q}\upharpoonright(\beta+1)\Vdash_{\P_{\beta+1}} \M\mbox{ is a weak }\lhd_\gamma\mbox{-chain w.r.t. }\G_{\beta+1}.$$

\item\claim There exist $r\in D$ and $s_0\in\P_{\beta+1}$ such that
\begin{parts}
    \item $\beta\in\dom(w_r)$,
    
    \item $\{N\in\M : N\cap\omega_1^V<M\cap\omega_1^V\}\subseteq\M_r$,
    
    \item $s_0\leq_{\P_{\beta+1}}{q}\upharpoonright(\beta+1),r\upharpoonright(\beta+1)$,
    
    \item $s_0\Vdash_{\P_{\beta+1}}r\in M^{\G_{\beta+1}}$.
\end{parts}
\begin{proof0}
    \item Let $\N:=\{N\in\M : N\cap\omega_1^V<M\cap\omega_1^V\}$ and let $G_{\beta+1}\leadsto V^{\P_{\beta+1}}$ arbitrary containing ${q}\upharpoonright(\beta+1)$.
    
    \item For all $N\in\N$, we have
    $$N\lhd_\gamma (M\downarrow\gamma)^{G_{\beta+1}}\cong_\gamma M^{G_{\beta+1}},$$
    which yields $P\in M^{G_{\beta+1}}$ such that $N=P\downarrow\gamma$; since $\gamma\in M^{G_{\beta+1}}$, we conclude $N\in M^{G_{\beta+1}}$.
    Hence, $\N\in M^{G_{\beta+1}}$.
    
    \item Let $D^*:=\{r\in D : r\upharpoonright(\beta+1)\in G_{\beta+1},\ \beta\in\dom(w_r),\ \N\subseteq\M_r\}\in M^{G_{\beta+1}}$.
    
    \item Note that $q\in D^*$.
    By elementarity, $M^{G_{\beta+1}}\models D^*\not=\emptyset$.
    Hence, there exists $r\in D\cap M^{G_{\beta+1}}$.
    
    \item By Forcing Theorem, there exists $s_0\in G_{\beta+1}$ such that
    \begin{parts}
        \item $s_0\leq q\upharpoonright(\beta+1),r\upharpoonright(\beta+1)$,
        
        \item $s_0\Vdash_{\P_{\beta+1}}r\in M^{\G_{\beta+1}}$.
    \end{parts}
    Conditions $r,s_0$ are as required.
\end{proof0}

\item\label{1662}\claim There exists $s\in\P_\gamma$ such that $s\leq_{\P_\gamma}q,r,s_0$.
\begin{proof0}
    
    \item Let $w_s:=w_{s_0}$, $\M_s:=\M_{s_0}\cup\M_q\cup\M_r$, $s:=(w_s,\M_s)$.
    We will be done once we show that $s\in\P_\gamma$.
    
    \item We have that $s\upharpoonright(\beta+1)=s_0\in\P_{\beta+1}$, so it remains to verify that
    \begin{eqnarray}
    s\upharpoonright(\beta+1) & \Vdash_{\P_{\beta+1}} & ({\M}_s^\gamma [\G_\beta]\mbox{ is a weak }\lhd_\gamma\mbox{-chain w.r.t. }\dot{G}_{\beta+1}), \label{1146}\\
    s\upharpoonright\beta & \Vdash_{\P_\beta} & (\forall N\in {\M}_s^\gamma[\G_\beta])(w_s(\beta)\mbox{ is semi-}N[\dot{G}_{\beta}]^{\dot{\Q}_\beta}\mbox{-generic}).\label{1147}
    \end{eqnarray}
    
    \item Let $G_{\beta+1}\leadsto V^{\P_\beta}$ be arbitrary containing $s\upharpoonright(\beta+1)$.
    Let us first verify that $\M_s^\gamma [G_\beta]$ is a weak $\lhd_\gamma$-chain w.r.t. $G_{\beta+1}$.
    
    \item Let $\M^*:=\{N\in\M : N\cap\omega_1^V\geq M\cap\omega_1^V\}$.
    Since $\{N\in\M_q^\gamma [G_\beta] : N\cap\omega_1^V<M\cap\omega_1^V\}\subseteq\M_r$, we have that 
    $$\M_s^\gamma [G_\beta]=\M_q^\gamma [G_\beta]\cup\M_r^\gamma [G_\beta]=\M_r^\gamma [G_\beta]\cup\M^*.$$
    
    \item Both $\M_r^\gamma [G_\beta]$ and $\M^*$ are weak $\lhd_\gamma$-chains w.r.t. $G_{\beta+1}$ and for all $N\in\M_r^\gamma [G_\beta]$ and all $P\in\M^*$, we have that
    $$N\cap\omega_1^V< M^{G_{\beta+1}}\cap\omega_1^V=M\cap\omega_1^V\leq P\cap\omega_1^V$$
    (``$<$" follows from $r\in M^{G_{\beta+1}}$ and ``$=$" from Claim \ref{1571}).
    Hence, it suffices to verify that for all $N\in\M_r^\gamma [G_\beta]$ and all $P\in\M^*$, we have that $N\lhd_\gamma P^{G_{\beta+1}}$.
    
    \item Since $r\in M^{G_{\beta+1}}$, we have $N\lhd_\gamma M^{G_{\beta+1}}$.
    
    \item Since $M\downarrow\gamma,P\in\M$ and $M\cap\omega_1^V\leq P\cap\omega_1^V$, we have
    $M\downarrow\gamma=P$ or $M\downarrow\gamma\lhd_\gamma P^{G_{\beta+1}}$.
    
    \item Hence,
    $$N\lhd_\gamma (M\downarrow\gamma)^{G_\beta+1}\subseteq P^{G_{\beta+1}},$$
    as required.
    
    \item Let us now verify that for every $N\in\M_s^\gamma [G_\beta]$, we have that $w_s(\beta)^{G_\beta}$ is semi-$N[G_\beta]^{\Q_\beta}$-generic.
    
    \item Either $N\in \M_q^\gamma [G_\beta]$ or $N\in \M_r^\gamma [G_\beta]$.
    
    \item If $N\in\M_q^\gamma [G_\beta]$, the conclusion follows from
    $$w_s(\beta)^{G_\beta}=w_{s_0}(\beta)^{G_\beta}\leq w_q(\beta)^{G_\beta};$$
    if $N\in\M_r^\gamma [G_\beta]$, the conclusion follows from
    $$w_s(\beta)^{G_\beta}=w_{s_0}(\beta)^{G_\beta}\leq w_r(\beta)^{G_\beta}.$$

    \item This concludes the proof of Claim \ref{1662}.
\end{proof0}

\item It is clear that $r$ and $s$ are as required in \ref{1645}.\qed 
\end{pfenum}

\subsubsection{Limit Case}
\label{limit}
\begin{pfenum}
\item We are assuming that $\gamma$ is a limit point of $\e$ and that
$$(\forall\alpha\in\e\cap\gamma)(\forall N\in\mathcal{C}_{>\alpha})(\forall q\in\P_\alpha)(N\downarrow\alpha\in\M_q\implies q\Vdash_{\P_\alpha}N^{\G_\alpha}\cap\omega_1^V=N\cap\omega_1^V).$$
Let us fix some $M\in\mathcal{C}_{>\gamma}$ and some $p\in\P_\gamma$ such that $M\downarrow\gamma\in\M_p$.
We want to show that $p$ is locally $M^{\P_\gamma}$-generic.

\item\label{1909} We use Local Genericity Criterion \ref{loc gen crit}.
Let $q\leq_{\P_\gamma}p$ be an arbitrary condition and let $D$ be an arbitrary dense open subset of $\P_\gamma$ such that
$$q\Vdash_{\P_\gamma}D\in M^{\G_\gamma}_{<\gamma}.$$
We want to find $r\in D$ and $s\leq_{\P_\gamma}q,r$ such that
$$s\Vdash_{\P_\gamma}\gamma\not\in M^{\G_\gamma}_{<\gamma}\vee r\in M^{\G_\gamma}_{<\gamma}.$$

\item\label{1563} By Lemma \ref{bounding activity}, we may assume, up to strengthening $q$, that there exist $\beta\in\e\cap\gamma$ and $\gamma^*\in (\beta,\gamma]$ such that
$$(\forall\xi\in\e\cap [\beta,\gamma))(q\upharpoonright\xi\Vdash_{\P_\xi}\sup(M^{\G_\xi}\cap\gamma)=\gamma^*).$$

\item If $q\Vdash_{\P_\gamma}\gamma\not\in M^{\G_\gamma}_{<\gamma}$, we can set $r=s$ to be an arbitrary element of $D$ which is $\leq q$.
Hence, we will assume that $q\not\Vdash_{\P_\gamma}\gamma\not\in M^{\G_\gamma}_{<\gamma}$.
Up to strengthening $q$, we may assume that $q\Vdash_{\P_\gamma}\gamma\in M^{\G_\gamma}_{<\gamma}$.
Up to further strengthening and increasing of $\beta$, we may assume that $q\Vdash_{\P_\gamma}\gamma\in M^{\G_\gamma}_{<\beta}$, or equivalently
$$q\upharpoonright\beta\Vdash_{\P_\beta}\gamma\in M^{\G_\beta}_{<\beta}.$$

\item Up to strengthening $q$ and increasing $\beta$, we may also assume that 
$$q\upharpoonright\beta\Vdash_{\P_\beta}D\in M^{\G_\beta}_{<\beta}.$$

\item\label{1572}\claim We may assume that for all $N\in\M_q$, either 
$${q}\upharpoonright\beta\Vdash_{\P_{\beta}}\sup(N^{\G_{\beta}}\cap\gamma^*)=\gamma^*$$
or
$$(\forall \xi\in\e\cap [\beta,\gamma^*))({q}\upharpoonright\xi\Vdash_{\P_\xi}\sup(N^{\G_\xi}\cap\gamma^*)\leq\beta).$$
\begin{proof0}
    \item Let $\M_q=\{N_i : i<k\}$, let $\eta_{-1}:=\beta$, and let $t_{-1}:=q\upharpoonright\gamma^*$.
    
    \item Let $i<k$.
    Suppose recursively that an ordinal $\eta_{i-1}\in\e\cap [\beta,\gamma^*)$ and a condition $t_{i-1}\leq_{\P_{\gamma^*}}q\upharpoonright\gamma^*$ satisfying 
    $$\M_{t_{i-1}}\cap\mathcal{C}_{\gamma^*}=\M_{q\upharpoonright\gamma^*}\cap\mathcal{C}_{\gamma^*}$$
    have been defined and that for all $j<i$, either
    $$t_{i-1}\upharpoonright\eta_{i-1}\Vdash_{\P_{\eta_{i-1}}}\sup(N_j^{\G_{\eta_{i-1}}}\cap\gamma^*)=\gamma^*$$
    or
    $$(\forall\xi\in\e\cap [\eta_{i-1},\gamma^*))(t_{i-1}\upharpoonright\xi\Vdash_{\P_\xi}\sup(N_j^{\G_\xi}\cap\gamma^*)\leq\eta_{i-1})).$$
    
    \item We want to define an ordinal $\eta_{i}\in\e\cap [\eta_{i-1},\gamma^*)$ and a condition $t_i\leq_{\P_{\gamma^*}}t_{i-1}$ satisfying
    $$\M_{t_i}\cap\mathcal{C}_{\gamma^*}=\M_{q\upharpoonright\gamma^*}\cap\mathcal{C}_{\gamma^*}$$
    such that for all $j\leq i$, either
    $$t_i\upharpoonright\eta_i\Vdash_{\P_{\eta_i}}\sup(N_j^{\G_{\eta_i}}\cap\gamma^*)=\gamma^*$$
    or
    $$(\forall\xi\in\e\cap [\eta_i,\gamma^*))(t_i\upharpoonright\xi\Vdash_{\P_\xi}\sup(N_j^{\G_\xi}\cap\gamma^*)\leq\eta_i)).$$
    
    \item If $N_i\in\mathcal{C}_{<\gamma^*}$, then let $\eta_i$ be the minimum of the set
    $$\e\cap [\max\{\eta_{i-1},\sup(N_i\cap\ord)\},\gamma^*)$$ and let $t_i:=t_{i-1}$.
    
    \item Otherwise, we have $N_i\in\mathcal{C}_{\geq\gamma^*}$.
    Then we can apply Lemma \ref{bounding activity} to
    \begin{parts}
        \item $\underaccent{\tilde}{\gamma}:=\gamma^*$,
        
        \item $\underaccent{\tilde}{M}:=N_i$,
        
        \item $\underaccent{\tilde}{p}:=t_{i-1}$,
    \end{parts}
    and obtain ordinals $\mu<\nu\leq\gamma^*$ and a condition $t_i\leq_{\P_{\gamma^*}}t_{i-1}$ satisfying
    $$\M_{t_i}\cap\mathcal{C}_{\gamma^*}=\M_{t_{i-1}}\cap\mathcal{C}_{\gamma^*},$$
    $$(\forall\xi\in\e\cap [\mu,\gamma^*))(t_i\upharpoonright\xi\Vdash_{\P_{\xi}}\sup(N_i^{\G_\xi}\cap\gamma^*)=\nu).$$
    
    \item If $\nu=\gamma^*$, then set
    $$\eta_i:=\max\{\eta_{i-1},\min(\e\cap [\mu,\gamma^*))\}.$$
    If $\nu<\gamma^*$, then set
    $$\eta_i:=\max\{\eta_{i-1},\min(\e\cap [\mu,\gamma^*)),\nu\}.$$
    
    \item By Lemma \ref{724}, we have $\nu\in\e$ and consequently $\eta_i\in\e$.
    
    \item This concludes the recursion.
    Note that $\M_{t_{k-1}}-\M_q\subseteq\mathcal{C}_{<\gamma^*}$, so there exists $\eta_k\in\e\cap [\eta_{k-1},\gamma^*)$ such that for all $N\in\M_{t_{k-1}}-\M_q$ and for all $\xi\in\e\cap [\eta_k,\gamma^*)$,
    $$\Vdash_{\P_\xi}\sup(N^{G_\xi}\cap\gamma^*)\leq\sup(N\cap\ord)\leq\eta_k.$$
    
    \item The claim now follows by replacing ${q}$ with $t_{k-1}q$ and by replacing $\beta$ with $\eta_{k}$.
\end{proof0}

\item Up to increasing $\beta$, we may assume that $\dom(w_q)\cap\gamma^*\subseteq\beta$.

\item For all $\xi\in\e\cap [\beta,\gamma^*)$, condition ${q}\upharpoonright\xi$ decides the value of $\M_q^{\gamma^*}[\G_\xi]$.
Up to increasing $\beta$, we may assume that all these values are equal to some fixed $\M\in [\mathcal{C}_{\gamma^*}]^{<\omega}$.
Hence, we can ensure the following:
\begin{parts}
    \item $(\forall\xi\in\e\cap [\beta,\gamma^*))({q}\upharpoonright\xi\Vdash_{\P_\xi}\M_q^{\gamma^*}[\G_\xi]=\M)$,
    
    \item ${q}\upharpoonright\beta\Vdash_{\P_{\beta}}(\M\mbox{ is a weak }\lhd_{\gamma^*}\mbox{-chain w.r.t. }\G_{\beta})$,
    
    \item $M\downarrow\gamma^*\in\M$.
\end{parts}

\item By \ref{1563} and Lemma \ref{724}, we have that
$${q}\upharpoonright\beta\Vdash_{\P_\beta}\sup(\e\cap M^{\G_\beta}\cap\gamma)=\gamma^*>\beta.$$
This means that we may strengthen $q$ and increase $\beta$ as to ensure
$${q}\upharpoonright\beta\Vdash_{\P_\beta}\beta\in M^{\G_\beta}_{<\beta}.$$

\item\claim There exist $r\in D$ and $s_0\in\P_\beta$ such that
\begin{parts}[ref=\arabic{pfenumi}$^\circ$\alph{partsi}]
    \item\label{1627} $(\forall N\in\M : N\cap\omega_1^V<M\cap\omega_1^V)(\exists N^+\in\M_r\cap\mathcal{C}_\gamma)(N=N^+\downarrow\gamma)$,
    
    \item $s_0\leq {q}\upharpoonright\beta,r\upharpoonright\beta$,
    
    \item $s_0\Vdash_{\P_\beta} r\in M^{\G_\beta}$.
\end{parts}
\begin{proof0}
    \item Let $G_\beta$ be a $V^{\P_\beta}$-generic containing ${q}\upharpoonright\beta$ and let us work inside $V[G_\beta]$.
    
    \item For $N\in\M$ with $N\cap\omega_1^V<M\cap\omega_1^V$, we may apply Lemma \ref{lifting vm} to
    \begin{parts}
        \item $\underaccent{\tilde}{M}:=M^{G_\beta}$,
        
        \item $\underaccent{\tilde}{\gamma}:=\gamma$,
        
        \item $\underaccent{\tilde}{N}:=N$,
    \end{parts}
    and obtain the unique $N^+\in M^{G_\beta}\cap\mathcal{C}_\gamma$ satisfying $N^+\downarrow\gamma^*=N$.
    Let
    $$\N:=\{N^+ : N\in\M,\, N\cap\omega_1^V<M\cap\omega_1^V\}\in M^{G_\beta}.$$
    
    \item Let $D^*:=\{r\in D : r\upharpoonright\beta\in G_\beta,\ \N\subseteq\M_r\}$.
    We want to find $r\in D^*\cap M^{G_\beta}$.
    
    \item Since $\beta\in M^{G_\beta}_{<\beta}$, model $M[G_\beta]$ is defined and $D^*\in M[G_\beta]$.
    By elementarity, it suffices to show that $D^*\not=\emptyset$.
    
    \item\textbf{Subclaim.} $\N$ is a weak $\lhd_\gamma$-chain w.r.t. $G_\beta$.
    \begin{proof0}
        \item Let $N,P\in\M$ be arbitrary satisfying $N\cap\omega_1^V\leq P\cap\omega_1^V<M\cap\omega_1^V$.
        We want to show that
        \begin{itemize}
            \item if $N\cap\omega_1^V= P\cap\omega_1^V$, then $N^+=P^+$,
            
            \item if $N\cap\omega_1^V< P\cap\omega_1^V$, then $N^+\lhd_\gamma (P^+)^{G_\beta}$.
        \end{itemize}
        The first case is obvious, since we have $N=P$.
        Let us consider the case when $N\cap\omega_1^V< P\cap\omega_1^V$.
        
        \item We have that $N\lhd_{\gamma^*}P^{G_\beta}$.
        Since $P\cong_{\gamma^*} P^+$, we have that $N\lhd_{\gamma^*}(P^+)^{G_\beta}$.
        
        \item Hence, there exists $N'\in (P^+)^{G_\beta}$ such that $N\cong_{\gamma^*}N'$.
        
        \item Since $P^+\in (M\downarrow\gamma)^{G_\beta}$, we conclude that 
        $$(P^+)^{G_\beta}\subseteq (M\downarrow\gamma)^{G_\beta}.$$
        
        \item This means that $N'\in (M\downarrow\gamma)^{G_\beta}$.
        Then there exists $N''\in M^{G_\beta}$ such that $N'\cong_\gamma N''$.
        
        \item Hence, $N'\downarrow\gamma=N''\downarrow\gamma\in M^{G_\beta}$, while $(N'\downarrow\gamma)\downarrow\gamma^*=N$.
        We conclude that 
        $$N^+=N'\downarrow\gamma\cong_{\gamma} N'\in (P^+)^{G_\beta},$$
        which means $N^+\lhd_\gamma(P^+)^{G_\beta}$.
    \end{proof0}
    

    
    \item There exists $t\in G_\beta$ such that
    \begin{parts}
        \item $t\leq {q}\upharpoonright\beta$,
        
        \item $t\Vdash_{\P_\beta}^V(\check{\N}$ is a weak $\lhd_\gamma$-chain w.r.t. $\G_\beta)$,
        
        \item $t\Vdash_{\P_\beta}^V(\forall P\in\check{\N})(P^{\G_\beta}\mbox{ is active at }\gamma^*)$.
    \end{parts}
    
    \item Let $w_u:=w_t$, let $\M_u:=\M_t\cup\N$, and let $u:=(w_u,\M_u)$.
    The previous point implies that $u\in\P_\gamma$.
    We also have that $u\upharpoonright\beta=t\in G_\beta$, which means that $u\in (\P_\gamma/\P_\beta)^{G_\beta}$.
    
    \item Since $\{v\in D : v\upharpoonright\beta\in G_\beta\}$ is dense in $(\P_\gamma/\P_\beta)^{G_\beta}$, we conclude that there exists $v\in D$ such that $v\upharpoonright\beta\in G_\beta$ and $v\leq_{\P_\gamma}u$; in particular, $\N\subseteq\M_v$ and consequently $v\in D^*$.
    
    \item Hence, there exists $r\in D^*\cap M[G_\beta]$.
    This means that
    \begin{parts}
        \item $r\in D$,
        
        \item $r\upharpoonright\beta\in G_\beta$,
        
        \item $\N\subseteq\M_r$.
    \end{parts}
    
    \item By Forcing Theorem, there exists $s_0\in G_\beta$ such that
    \begin{parts}
        \item $s_0\leq q\upharpoonright\beta,r\upharpoonright\beta$,
        
        \item $s_0\Vdash_{\P_\beta}r\in M^{\G_\beta}$.
    \end{parts}
\end{proof0}

\item\label{2098}\claim There exists $s\in\P_\gamma$ such that $s\leq_{\P_\gamma}q,r,s_0$.
\begin{proof0}
    \item Let $t:=s_0 q\in \P_\gamma$.
    Note that
    $$\dom(w_t)\cap [\beta,\gamma^*)=\dom(w_q)\cap [\beta,\gamma^*)=\emptyset.$$

    \item Let $\dom(w_s):=\dom(w_t)\cup\dom(w_r)$.
    For $\xi\in\dom(w_s)- [\beta,\gamma^*)$, let $w_s(\xi):=w_t(\xi)$.

    \item For $\xi\in\dom(w_s)\cap [\beta,\gamma^*)=\dom(w_r)-\beta$, let us denote by $\eta$ its successor in $\e$.
    
    \item Since $s_0\Vdash_{\P_\beta}r\in M^{\G_\beta}$ and $t\upharpoonright\xi\leq s_0$, we have that $t\upharpoonright\xi\Vdash_{\P_\xi}(M^{\G_\xi}\mbox{ is active at }\eta)$.
    
    \item Hence, we may apply Lemma \ref{1372} to 
    \begin{parts}
        \item $\underaccent{\tilde}{\mu}:=\xi$, $\underaccent{\tilde}{\nu}:=\eta$,
        
        \item $\underaccent{\tilde}{p}:=t\upharpoonright\eta$,
        
        \item $\underaccent{\tilde}{M}:=M\downarrow\eta$,
        
        \item $\underaccent{\tilde}{\dot{u}}:=w_r(\xi)$
    \end{parts}
    and obtain a canonical $\P_\xi$-name $w_s(\xi)$ for an element of $\dot{\Q}_\xi$ such that
    $$t\upharpoonright\xi\Vdash_{\P_\xi}w_s(\xi)\leq w_r(\xi),$$
    $$t\upharpoonright\xi\Vdash_{\P_\xi}(\forall N\in\M_t^\eta [\G_\xi])(N\cap\omega_1^V\geq M\cap\omega_1^V\implies (w_s(\xi)\mbox{ is semi-}N[\G_\xi]^{\dot{\Q}_\xi}\mbox{-generic})).$$
    
    \item Let $\M_s:=\M_t\cup\M_r$ and let $s:=(w_s,\M_s)$.
    We will be done if we show $s\in\P_\gamma$.
    We show by induction on $\eta\in\e\cap [\beta,\gamma]$ that $s\upharpoonright\eta\in\P_\eta$.
    
    \item\textbf{Case.} {$\eta=\beta$}
    \begin{proof}
    This case is obvious since $s\upharpoonright\beta=s_0$.
    \end{proof}

    \item\label{1820}\textbf{Case.} Let us assume that
    {$\eta\in (\beta,\gamma^*]$} and that {$\eta$ has a predecessor $\xi$ inside $\e$}.
    \begin{proof0}
        \item We need to verify that
        $$\xi\in\dom(w_s)\implies s\upharpoonright\xi\Vdash_{\P_\xi}(\forall N\in\M_s^\eta [\G_\xi])(w_s(\xi)\mbox{ is semi-}N[\G_\xi]^{\dot{\Q}_\xi}\mbox{-generic})$$
        and that
        $$s\upharpoonright(\xi+1)\Vdash_{\P_{\xi+1}}(\M_s^\eta [\G_\xi]\mbox{ is a weak }\lhd_\eta\mbox{-chain w.r.t. }\G_{\xi+1}).$$
        
        \item Let $G_\xi\leadsto V^{\P_\xi}$ be an arbitrary generic containing $s\upharpoonright\xi$ and let $\Q_\xi:=\dot{\Q}_\xi^{G_\xi}$.
        Note that 
        $$\M_s^\eta [G_\xi]=\M_t^\eta [G_\xi]\cup\M_r^\eta [G_\xi].$$
        
        \item Since $\M_t^\eta [G_\xi]\subseteq\M_t\cap\mathcal{C}_{\geq\eta}=\M_q$, we have that $\M_t^\eta [G_\xi]=\M_q^\eta [G_\xi]$ and consequently
        $$\M_s^\eta [G_\xi]=\M_q^\eta [G_\xi]\cup\M_r^\eta [G_\xi].$$
        
        \item\label{1755}\textbf{Subclaim.} $\M_q^\eta [G_\xi]\subseteq\M\downarrow\eta$
        \begin{proof0}
            \item Let $N\in\M_q^\eta [G_\xi]$ be arbitrary.
            Then there exists $P\in\M_q\cap\mathcal{C}_{\geq\eta}$ such that $\xi\in P^{G_\xi}$ and $N=P\downarrow\eta$.
            
            \item We now have $\sup(P^{G_\xi}\cap\gamma)>\xi\geq\beta$ and it follows from Claim \ref{1572} that 
            $$\sup(P^{G_\xi}\cap\gamma)=\gamma^*.$$
            
            \item This means that $P\downarrow\gamma^*\in\M_q^{\gamma^*}[G_\xi]=\M$, while $N=(P\downarrow\gamma^*)\downarrow\eta$.
        \end{proof0}
        
        \item\label{1781}\textbf{Subclaim.} $\{N\in\M_q^\eta [G_\xi] : N\cap\omega_1^V< M\cap\omega_1^V\}\subseteq\M_r^\eta [G_\xi]$
        \begin{proof0}
            \item Let $N\in\M_q^\eta [G_\xi]$ be such that $N\cap\omega_1^V< M\cap\omega_1^V$.
            
            \item By Subclaim \ref{1755}, there exists $P\in\M$ such that $N=P\downarrow\eta$.
            
            \item We have $P\cap\omega_1^V< M\cap\omega_1^V$, so Claim \ref{1627} implies that there exists $Q\in\M_r\cap\mathcal{C}_\gamma$ such that $P=Q\downarrow\gamma^*$.
            
            \item Note that
            $$\xi\in N^{G_\xi}\cong_\eta P^{G_\xi}\cong_{\gamma^*}Q^{G_\xi},$$
            so $Q\downarrow\eta\in\M_r^\eta [G_\xi]$, while $N=Q\downarrow\eta$.
        \end{proof0}
        
        \item\textbf{Subclaim.} For all {$\xi\in\dom(w_s)$} and for all {$N\in\M_s^\eta [G_\xi]$}, we have that 
        $w_s(\xi)^{G_\xi}$ is semi-$N[G_\xi]^{\Q_\xi}$-generic.
        \begin{proof0}
            \item Suppose first that $N\in\M_r^\eta [G_\xi]$.
            Since $\xi\in [\beta,\gamma^*)$, we have $\xi\in\dom(w_r)-\beta$ and consequently that $w_r(\xi)^{G_\xi}$ is semi-$N[G_\xi]^{\Q_\xi}$-generic.
            The conclusion now follows from the fact that $w_s(\xi)^{G_\xi}\leq w_r(\xi)^{G_\xi}$.
            
            \item Suppose next that $N\in\M_q^\eta [G_\xi]$ and that $N\cap\omega_1^V\geq M\cap\omega_1^V$.
            Then $w_s(\xi)^{G_\xi}$ is semi-$N[G_\xi]^{\Q_\xi}$-generic by the choice of $w_s(\xi)$.
            
            \item Suppose finally that $N\in\M_q^\eta [G_\xi]$ and that $N\cap\omega_1^V< M\cap\omega_1^V$.
            By Subclaim \ref{1781}, we have that $N\in\M_r^\eta [G_\xi]$, which means that $w_r(\xi)^{G_\xi}$ is semi-$N[G_\xi]^{\Q_\xi}$-generic.
            The desired conclusion now follows from the fact that $w_s(\xi)^{G_\xi}\leq w_r(\xi)^{G_\xi}$.
        \end{proof0}
        
        \item\textbf{Subclaim.} Let $G_{\xi+1}\leadsto V^{\P_{\xi+1}}$ be an arbitrary generic extending $G_\xi$ and containing $s\upharpoonright(\xi+1)$.
        Then $\M_s^\eta [G_\xi]$ is a weak $\lhd_\eta$-chain w.r.t. $G_{\xi+1}$.
        \begin{proof0}
            \item Let $\M^*:=\{P\in\M_q^\eta [G_\xi] : P\cap\omega_1^V\geq M\cap\omega_1^V\}$.
            By Subclaim \ref{1781}, set $\M_s^\eta [G_\xi]$ is the union of sets $\M_r^\eta [G_\xi]$ and $\M^*$.
            
            \item Each of sets $\M_r^\eta [G_\xi]$ and $\M^*$ is a weak $\lhd_\eta$-chain w.r.t. $G_{\xi+1}$.
            Furthermore, for every $N\in \M_r^\eta [G_\xi]$ and every $P\in\M^*$, we have
            $$N\cap\omega_1^V<M\cap\omega_1^V\leq P\cap\omega_1^V.$$
            
            \item Hence, it remains to verify that for every $N\in \M_r^\eta [G_\xi]$ and every $P\in\M^*$, we have $N\lhd_\eta P^{G_{\xi+1}}$.
            
            \item Since $N\in \M_r^\eta [G_\xi]$ and $r\in M^{G_\xi}$, we have $N\lhd_\eta  M^{G_\xi}$.
            
            \item By Subclaim \ref{1755}, there exists $Q\in\M$ such that $P=Q\downarrow\eta$.
            Since $Q\cap\omega_1^V\geq M\cap\omega_1^V$, we conclude that $M\downarrow\gamma^*=Q$ or $M\downarrow\gamma^*\lhd_{\gamma^*}Q^{G_\xi}$.
            
            \item Hence, there exists $R\in Q^{G_\xi}\cup \{Q\}$ such that $M\cong_{\gamma^*}R$.
            
            \item We conclude $N\lhd_\eta R^{G_\xi}\subseteq Q^{G_\xi}$.
            
            \item Thus,
            $$N\lhd_\eta Q^{G_\xi}\cong_\eta P^{G_\xi}\subseteq P^{G_{\xi+1}},$$
            as required.
        \end{proof0}
        
        \item This concludes the proof of Case \ref{1820}
    \end{proof0}
    
    \item\label{1823}\textbf{Case.} Let us assume that $\eta\in (\beta,\gamma^*]$ and that $\eta$ is a limit point of $\e$.
    \begin{proof0}
        \item We need to show that there exists $\beta_0<\eta$ such that
        $$(\forall\xi\in\e\cap(\beta_0,\eta))(s\upharpoonright(\xi+1)\Vdash_{\P_{\xi+1}}(\M_s^\eta [\G_\xi]\mbox{ is a weak }\lhd_\eta\mbox{-chain w.r.t. }\G_{\xi+1})).$$
        We will show that $\beta_0:=\beta$ works.
        
        \item Let $\xi\in\e\cap(\beta,\eta)$ and let $G_\xi\leadsto V^{\P_\xi}$ satisfying $s\upharpoonright\xi\in G_\xi$.
        
        \item We have $\M_s^\eta [G_\xi]=\M_t^\eta [G_\xi]\cup\M_r^\eta [G_\xi]=\M_q^\eta [G_\xi]\cup\M_r^\eta [G_\xi]$.
        
        \item\label{1833}\textbf{Subclaim.} $\M_q^\eta [G_\xi]\subseteq\M\downarrow\eta$
        \begin{proof0}
            \item Let $N\in\M_q^\eta [G_\xi]$ be arbitrary.
            Then there exists $P\in\M_q\cap\mathcal{C}_{\geq\eta}$ such that 
            $$\sup(P^{G_\xi}\cap\eta)=\eta$$
            and $N=P\downarrow\eta$.
            
            \item We now have $\sup(P^{G_\xi}\cap\gamma)>\beta$ and it follows from Claim \ref{1572} that 
            $$\sup(P^{G_\xi}\cap\gamma)=\gamma^*.$$
            
            \item This means that $P\downarrow\gamma^*\in\M_q^{\gamma^*}[G_\xi]=\M$, while $N=(P\downarrow\gamma^*)\downarrow\eta$.
        \end{proof0}
        
        \item\label{1846}\textbf{Subclaim.} $\{N\in\M_q^\eta [G_\xi] : N\cap\omega_1^V< M\cap\omega_1^V\}\subseteq\M_r^\eta [G_\xi]$
        \begin{proof0}
            \item Let $N\in\M_q^\eta [G_\xi]$ be such that $N\cap\omega_1^V< M\cap\omega_1^V$.
            
            \item By Subclaim \ref{1833}, there exists $P\in\M$ such that $N=P\downarrow\eta$.
            
            \item We have $P\cap\omega_1^V< M\cap\omega_1^V$, so Claim \ref{1627} implies that there exists $Q\in\M_r\cap\mathcal{C}_\gamma$ such that $P=Q\downarrow\gamma^*$.
            
            \item Since
            $$N^{G_\xi}\cong_\eta P^{G_\xi}\cong_{\gamma^*}Q^{G_\xi},$$
            we have that 
            $$\sup(Q^{G_\xi}\cap\eta)=\sup(N^{G_\xi}\cap\eta)=\eta,$$
            i.e. $Q\downarrow\eta\in\M_r^\eta [G_\xi]$.
            
            \item The conclusion now follows from the fact that $N=Q\downarrow\eta$.
        \end{proof0}
        
        \item\textbf{Subclaim.} Let $G_{\xi+1}\leadsto V^{\P_{\xi+1}}$ be an arbitrary generic extending $G_\xi$ and containing $s\upharpoonright(\xi+1)$.
        Then $\M_s^\eta [G_\xi]$ is a weak $\lhd_\eta$-chain w.r.t. $G_{\xi+1}$.
        \begin{proof0}
            \item Let $\M^*:=\{P\in\M_q^\eta [G_\xi] : P\cap\omega_1^V\geq M\cap\omega_1^V\}$.
            By Subclaim \ref{1846}, set $\M_s^\eta [G_\xi]$ is the union of sets $\M_r^\eta [G_\xi]$ and $\M^*$.
            
            \item Each of sets $\M_r^\eta [G_\xi]$ and $\M^*$ is a weak $\lhd_\eta$-chain w.r.t. $G_{\xi+1}$.
            Furthermore, for every $N\in \M_r^\eta [G_\xi]$ and every $P\in\M^*$, we have
            $$N\cap\omega_1^V<M\cap\omega_1^V\leq P\cap\omega_1^V.$$
            
            \item Hence, it remains to verify that for every $N\in \M_r^\eta [G_\xi]$ and every $P\in\M^*$, we have $N\lhd_\eta P^{G_{\xi+1}}$.
            
            \item Since $N\in \M_r^\eta [G_\xi]$ and $r\in M^{G_\xi}$, we have $N\lhd_\eta  M^{G_\xi}$.
            
            \item By Subclaim \ref{1833}, there exists $Q\in\M$ such that $P=Q\downarrow\eta$.
            Since $Q\cap\omega_1^V\geq M\cap\omega_1^V$, we conclude that $M\downarrow\gamma^*=Q$ or $M\downarrow\gamma^*\lhd_{\gamma^*}Q^{G_\xi}$.
            
            \item Hence, there exists $R\in Q^{G_\xi}\cup \{Q\}$ such that $M\cong_{\gamma^*}R$.
            
            \item We conclude $N\lhd_\eta R^{G_\xi}\subseteq Q^{G_\xi}$.
            
            \item Thus,
            $$N\lhd_\eta Q^{G_\xi}\cong_\eta P^{G_\xi}\subseteq P^{G_{\xi+1}},$$
            as required.
        \end{proof0}
    
        \item This concludes the proof of Case \ref{1823}.
    \end{proof0}
    
    \item\textbf{Case.} Let us assume that $\eta\in (\gamma^*,\gamma]$.
    \begin{proof0}
        \item\label{2064} Since $s_0\Vdash_{\P_\beta}r\in M^{\G_\beta}\cap\P_\gamma$, we have that $\dom(w_r)\subseteq\gamma^*$ and $\M_r\subseteq\mathcal{C}_{<\gamma^*}$.
        Hence,
        $$w_{s\upharpoonright\eta}=(w_s\upharpoonright\gamma^*)\cup(w_t\upharpoonright [\gamma^*,\eta))=(w_s\upharpoonright\gamma^*)\cup(w_q\upharpoonright [\gamma^*,\eta)),$$
        $$\M_{s\upharpoonright\eta}=(\M_s\cap\mathcal{C}_{\leq\gamma^*})\cup (\M_t\downarrow\eta)=(\M_s\cap\mathcal{C}_{\leq\gamma^*})\cup (\M_q\downarrow\eta).$$
        
        \item Suppose first that $\eta$ is the successor of some $\beta$ inside $\e$.
        We have to show that
        \begin{equation} 
        \label{2072}
            \xi\in\dom(w_s)\implies s\upharpoonright\xi\Vdash_{\P_\xi}(\forall N\in\M_s^\eta [\G_\xi])(w_s(\xi)\mbox{ is semi-}N[\G_\xi]^{\dot{\Q}_\xi}\mbox{-generic})
        \end{equation}
        and that
        \begin{equation}
        \label{2077}
            s\upharpoonright(\xi+1)\Vdash_{\P_{\xi+1}}(\M_s^\eta [\G_\xi]\mbox{ is a weak }\lhd_\eta\mbox{-chain w.r.t. }\G_{\xi+1}).
        \end{equation}
        
        \item By \ref{2064}, we have that $\Vdash_{\P_\xi}\M_s^\eta [\G_\xi]=\M_q^\eta [\G_\xi]$ and $w_s(\xi)=w_q(\xi)$.
        Properties (\ref{2072}) and (\ref{2077}) now follow from the fact that $q$ is a condition and that $s\upharpoonright(\xi+1)\leq_{\P_{\xi+1}}q\upharpoonright(\xi+1)$.
        
        \item Suppose now that $\eta$ is a limit point of $\e$.
        We have to show that there exists $\bar{\eta}<\eta$ such that for all $\xi\in\e\cap (\bar{\eta},\eta)$,
        \begin{equation}
        \label{2087}
            s\upharpoonright(\xi+1)\Vdash_{\P_{\xi+1}}(\M_s^\eta [\G_\xi]\mbox{ is a weak }\lhd_\eta\mbox{-chain w.r.t. }\G_{\xi+1}).
        \end{equation}
        
        \item We claim that $\bar{\eta}:=\gamma^*$ works.
        Namely, for all $\xi\in\e\cap(\gamma^*,\eta)$, we have that 
        $$\Vdash_{\P_\xi}\M_s^\eta [\G_\xi]=\M_q^\eta [\G_\xi],$$
        $$s\upharpoonright(\xi+1)\leq_{\P_{\xi+1}}q\upharpoonright(\xi+1),$$
        so the conclusion again follows from the fact that $q$ is a condition.
    \end{proof0}
    
    \item This concludes the proof of Claim \ref{2098}.
\end{proof0}

\item It is now clear that $s$ and $r$ as required in \ref{1909}, concluding the argument.\qed
\end{pfenum}

\subsection{Semi-properness and Chain Condition}
\begin{enumerate}
\item\label{semi-genericity}\lemma Let $\gamma\in\e$, let $M\in\mathcal{C}_{>\gamma}$, and let $p\in\P_\gamma$ satisfy $M\downarrow\gamma\in\M_p$.
Then 
$$p\Vdash_{\P_\gamma}M^{\G_\gamma}\cap\omega_1^V=M\cap\omega_1^V.$$
\begin{proof0}
    \item This is shown by induction on $\gamma$.
    
    \item\textbf{Case.} $\gamma=\min\e$
    \begin{proof}
    By Transfer Theorem, condition $p$ is locally $M^{\P_\gamma}$-generic, which in this case translates to $M^{\P_\gamma}$-generic.
    \end{proof}
    
    \item\textbf{Case.} Suppose that $\gamma$ is a successor of some $\beta$ inside $\e$.
    \begin{proof0}
        \item Let $G_\gamma\leadsto V^{\P_\gamma}$ be arbitrary containing $p$ and let us verify that 
        $$M^{G_\gamma}\cap\omega_1^V=M\cap\omega_1^V.$$
        
        \item By the IH and Transfer Theorem, we have that
        $$M^{G_\beta}\cap\omega_1^V=M\cap\omega_1^V,\quad M^{G_\gamma}=M^{G_{\beta+1}}.$$
        Hence, it suffices to verify that
        $$M^{G_{\beta+1}}\cap\omega_1^V=M^{G_\beta}\cap\omega_1^V.$$
        
        \item If $\beta\not\in M^{G_{\beta}}_{<\beta}$, then
        $$M^{G_\beta}_{<\beta}=M^{G_\beta}=M^{G_\beta+1}.$$
        Let us then assume that $\beta\in M^{G_{\beta}}_{<\beta}$.
        
        \item Since conditions $q\in\P_\gamma$ such that $\beta\in\dom(w_q)$ are dense below $p$, there exists $q\in G_\gamma$ such that $q\leq p$ and $\beta\in\dom(w_q)$.
        
        \item We have $M\downarrow\gamma\in\M_q^\gamma [G_\beta]$.
        Since $\beta\in\dom(w_q)$, we conclude that $w_q(\beta)^{G_\beta}$ is semi-$(M\downarrow\gamma)[G_\beta]^{\Q_{\beta}}$-generic.
        
        \item Since $q\in G_\beta$, we conclude that
        $$(M\downarrow\gamma)[G_{\beta+1}]\cap\omega_1^{V[G_\beta]}=(M\downarrow\gamma)[G_\beta]\cap\omega_1^{V[G_\beta]}.$$
        This suffices for the conclusion.
    \end{proof0}
    
    \item\textbf{Case.} $\gamma$ is a limit point of $\e$.
    \begin{proof}
    By the IH and Transfer Theorem, we have that 
    $$p\Vdash_{\P_\gamma}M^{\G_\gamma}=M^{\G_\gamma}_{<\gamma}.$$
    Hence,
    $$p\Vdash_{\P_\gamma}M^{\G_\gamma}\cap\omega_1^V=\bigcup_{\alpha\in\e\cap\gamma}(M^{G_\alpha}\cap\omega_1^V)=M\cap\omega_1^V$$
    by another application of the IH.
    \end{proof}
    
    \item This concludes the induction.
\end{proof0}

\item\proposition For all $\gamma\in\e^*$, poset $\P_\gamma$ is semi-proper.
\begin{proof0}
    \item Suppose first that $\gamma\in\e$.
    It suffices to show that for every $M\in\mathcal{C}_{>\gamma}$ and every $p\in M$, there exists $q\leq p$ which is semi-$M^{\P_\gamma}$-generic.
    
    \item By Lemma \ref{adding a side condition}, there exists $q\leq p$ such that $M\downarrow\gamma\in\M_q$.
    By the previous lemma, condition $q$ is semi-$M^{\P_\gamma}$-generic.
    
    \item Suppose now that $\gamma\in\e^+$.
    Then
    $$\P_\gamma\simeq \P_{\gamma-1}*\dot{\Q}_{\gamma-1},$$
    which suffices to conclude semi-properness.
\end{proof0}

\item\lemma Suppose that
\begin{parts}[label=\arabic*.]
    \item $\gamma\in\e$,

    \item $\cof(\gamma)=\omega_1$,

    \item $p\in\P_\gamma$.
\end{parts}
Then there exist $q$ and $\bar\gamma$ such that
\begin{parts}
    \item $q\leq_{\P_\gamma}p$,
    
    \item $\bar\gamma\in\e\cap\gamma$,

    \item $\dom(w_q)\subseteq\bar\gamma$,

    \item for all $\alpha\in\e\cap [0,\gamma]$, for all $N\in\M_{q}$, $q\upharpoonright\alpha\Vdash_{\P_\alpha}\sup(N^{\G_\alpha}\cap\gamma)\leq\bar\gamma$.


\end{parts}
\begin{proof0}
    \item We note that $\gamma$ is necessarily a limit point of $\e$.
    
    \item Let $M\in\mathcal{C}_{>\gamma}$ be such that $\gamma,p\in M$ and let $q\leq p$ be such that $\dom(w_{q})=\dom(w_p)$ and $\M_{q}=\M_p\cup\{M\downarrow\gamma\}$ (cf. Lemma \ref{adding a side condition}).
    
    \item Let $\bar{\gamma}:=\sup(\gamma\cap M)<\gamma$.
    We have that $\cof(\bar{\gamma})=\omega$, that $\bar{\gamma}$ is a limit point of $\e$, and that $\dom(w_{q})=\dom(w_p)\subseteq\bar{\gamma}$.
    
    \item\claim $q\Vdash_{\P_\gamma}\sup(M^{\G_\gamma}\cap\gamma)=\bar{\gamma}$
    \begin{proof0}
        \item Let $G_\gamma\leadsto V^{\P_\gamma}$ be an arbitrary generic containing $q$ and let us show that 
        $$\sup(M^{G_\gamma}\cap\gamma)=\bar{\gamma}.$$
        
        \item Let $f\in M$ map cofinally $\omega_1^V$ into $\gamma$.
        
        \item By Lemma \ref{semi-genericity}, we have that $M^{G_\gamma}\cap\omega_1^V=M\cap\omega_1^V$.
        
        \item Hence,
        $$\sup(M^{G_\gamma}\cap\gamma)=\sup(f[M^{G_\gamma}\cap\omega_1^V])=\sup(f[M\cap\omega_1^V])=\sup(M\cap\gamma)=\bar{\gamma}.$$
    \end{proof0}
    
    \item\label{2199} Consequently, for all $\alpha\in\e\cap [0,\gamma]$ and for all $N\in\M_{q}$, we have
    $$q\upharpoonright\alpha\Vdash_{\P_\alpha}\sup(N^{\G_\alpha}\cap\gamma)\leq\bar{\gamma}.$$
    
        
        
        
    
\end{proof0}

\item\theorem $\P_\kappa$ has $\kappa$-c.c.
\begin{proof0}
    \item Assume otherwise.
    
    \item Let $S:=\{\alpha\in\e : \cof(\alpha)=\omega_1\}$.
    Then there exists an anti-chain $(p_\alpha : \alpha\in S)$ in $\P_\kappa$.
    
    \item By the previous lemma, for all $\alpha\in S$, there exist $q_\alpha$ and $\bar\alpha$ such that
    \begin{parts}
        \item $q_\alpha\leq_{\P_\alpha}p_\alpha\upharpoonright\alpha$,

        \item $\bar\alpha\in\e\cap\alpha$,

        \item $\dom(w_{q_\alpha})\subseteq\bar\alpha$,

        \item for all $\xi\in\e\cap [0,\alpha]$, for all $N\in\M_{q_\alpha}$, $q_\alpha\upharpoonright\xi\Vdash_{\P_\xi}\sup(N^{\G_\xi}\cap\alpha)\leq\bar\alpha$.


    \end{parts}
    
    \item By pressing down, there exist a stationary $S_1\subseteq S$ and $\gamma\in\e$ such that for all $\alpha\in S_1$, we have $\bar\alpha=\gamma$.
    
    
    \item There exist $\alpha,\beta\in S_1$ such that
    \begin{parts}
        \item $\gamma<\alpha<\beta$,

        \item $q_\alpha\in\P_\beta$,
        
        \item $q_\alpha\upharpoonright\gamma=q_\beta\upharpoonright\gamma$.
    \end{parts}
    
    
    \item Let $q:=q_\alpha\upharpoonright\gamma=q_\beta\upharpoonright\gamma$ and let
    \begin{eqnarray*}
    w_r & := & w_q\cup(w_{p_\alpha}\upharpoonright [\gamma,\beta))\cup(w_{p_\beta}\upharpoonright [\beta,\kappa)),\\
    \M_r & := & \M_q\cup \M_{p_\alpha}\cup \M_{p_\beta},\\
    r & := & (w_r,\M_r).
    \end{eqnarray*}
    It is now easily seen that $r\in\P_\kappa$ and $r\leq p_\alpha,\, p_\beta$, which contradicts $p_\alpha\bot_{\P_\kappa} p_\beta$.
\end{proof0}

\item\theorem $\P_\kappa$ is semi-proper.
\begin{proof0}
    \item Let $\theta$ be large enough regular, let $M\prec(H_\theta,\in,\kappa,U)$ be countable, let $p\in\P_\kappa\cap M$, and let us find $q\leq_{\P_\kappa} p$ which is semi-$M^{\P_\kappa}$-generic.
    
    \item Let $\bar{\kappa}:=\sup(\kappa\cap M)\in\e\cap\kappa$, let $N:=M\cap V_\kappa\in\mathcal{C}_{\bar{\kappa}}$, and let $q\in\P_{\bar{\kappa}}$ be such that
    \begin{parts}
        \item $q\leq p$,
        
        \item $\dom(w_q)=\dom(w_p)$,
        
        \item $\M_q=\M_p\cup\{N\}$
    \end{parts}
    (cf. Lemma \ref{adding a side condition}; note that $p\in\P_{\bar{\kappa}}$).
    We want to show that $q$ is semi-$M^{\P_{\kappa}}$-generic.
    
    \item Let $\tau\in M^{\P_\kappa}$ with $\Vdash_{\P_\kappa}\tau<\omega_1^V$ be arbitrary.
    We want to show that $q\Vdash_{\P_{\kappa}}\tau\in M$.
    
    \item There exists $A\in M$ such that
    \begin{parts}
        \item $A$ is a maximal anti-chain of $\P_\kappa$,
        
        \item for all $r\in A$, there exists unique $\alpha_r<\omega_1^V$ such that $r\Vdash_{\P_\kappa}\tau=\alpha_r$.
    \end{parts}
    
    \item By the previous theorem, there exists $\gamma\in\e\cap M$ such that $(w_q,\M_p)\in\P_\gamma$ and $A\subseteq\P_\gamma$.
    
    \item We can now find $\sigma\in N^{\P_{\gamma}}$ such that for all $r\in A$, we have $r\Vdash_{\P_{\gamma}}\sigma=\alpha_r$.
    Note that then $\Vdash_{\P_\kappa}\tau=\sigma$ and consequently $\Vdash_{\P_{\gamma}}\sigma<\omega_1^V$.
    
    \item Let $q_\tau:=(w_q,\M_p\cup\{N\downarrow\gamma\})$.
    It is clear that $q_\tau\in\P_\gamma$ and $q\leq_{\P_\kappa}q_\tau$.
    
    \item By Lemma \ref{semi-genericity}, we have that $q_\tau$ is semi-$N^{\P_\gamma}$-generic.
    In particular, $q_\tau\Vdash_{\P_\gamma}\sigma\in N$.
    
    \item Hence,
    $$q\leq_{\P_\kappa}q_\tau\Vdash_{\P_\kappa}\tau=\sigma\in N\subseteq M,$$
    as required.
\end{proof0}

\end{enumerate}
\newpage

\section{Saturating \texorpdfstring{$\ns$}{NS}}

\subsection{Careful Collapse}
\begin{enumerate}

\item\notation For $M$ a virtual model, we denote $\delta_M:=\min(\ord-M)$.

\item\definition 
    Suppose that
    \begin{assume}
        \item $\theta>\omega$ is regular,

        \item $U$ is a stationary subset of $[H_\theta]^\omega$.
    \end{assume}
    \textit{The collapsing poset $\col_U$ guided by $U$}, is defined as follows.
    \begin{parts}
        \item $p\in\col_U$ if and only if $p=(\M_p,d_p)$ where
        \begin{parts}
            \item for all $M\in\M_p$, either $M\in U$ or there does not exist $N\in U$ satisfying $M\subseteq N$ and $\delta_N=\delta_M$,
            
            \item $\M_p$ is a finite $\in$-chain,

            \item $d_p:\M_p\to [H_\theta]^{<\omega}$,

            \item for all $M,N\in\M_p$, if $M\in N$, then $d_p(M)\subseteq N$.
        \end{parts}

        \item $p\leq q$ in $\col_U$ if and only if $\M_p\supseteq\M_q$ and for all $M\in\M_q$, $d_p(M)\supseteq d_q(M)$.
    \end{parts}

\item\lemma
    \label{90''}
    Suppose that
    \begin{assume}
        \item $\theta>\omega$ is regular,

        \item $U$ is a stationary subset of $[H_\theta]^\omega$,

        \item $x\in H_\theta$,

        \item $D:=\{p\in\col_U : \exists M\in\M_p,x\in M\}$.
    \end{assume}
    Then $D$ is dense in $\col_U$.
\begin{proof0}
    \item Let $p_0\in \col_U$ be arbitrary and let us find $p\in D$ satisfying $p\leq p_0$.

    \item Since $U$ is stationary, there exists $M\in U$ such that $x,p_0\in M$.

    \item Then $p:=(\M_{p_0}\cup\{M\},d_{p_0}\cup\{(M,\emptyset)\})$ is as required.
\end{proof0}

\item\lemma
    Suppose that
    \begin{assume}
        \item $\theta>\omega$ is regular,

        \item $U$ is a stationary subset of $[H_\theta]^\omega$,

        \item $G\leadsto V^{\col_U}$,

        \item $\M_G:=\cup\{\M_p : p\in G\}$.
    \end{assume}
    Then 
    \begin{parts}
        \item\label{98''} $\M$ is an $\in$-chain,

        \item\label{113''} $\otp(\M,\in)\leq\omega_1^V$,

        \item\label{115''} $\cup\M=H_{\theta}$.
    \end{parts}
\begin{proof0}
    \item Part \ref{98''}. is obvious.
    
    \item Let $(M_\xi : \xi<\alpha)$ be the $\in$-increasing enumeration of $\M$.
    The following mapping is strictly increasing:
    $$\alpha\to\omega_1^V:\xi\mapsto \delta_{M_\xi}.$$
    This shows part \ref{113''}.

    \item Lemma \ref{90''}. immediately implies part \ref{115''}.
\end{proof0}

\item\definition
    Suppose that
    \begin{assume}
        \item $\theta>\omega$ is regular,

        \item $U$ is a stationary subset of $[H_\theta]^\omega$,

        \item $g\leadsto V^{\col_U}$.
    \end{assume}
    Then we define
    \begin{parts}
        \item $\M_g:=\cup\{\M_p : p\in g\}$,

        \item $\alpha_g:=\otp(\M_g,\in)$,

        \item for all $\xi<\alpha_g$, $\M_g(\xi)$ is the $\xi^\mathrm{th}$ element of $(\M_g,\in)$.
    \end{parts}

\item\proposition
\label{continuity''}
    Suppose that 
    \begin{assume}
        \item $\theta>\omega$ is regular,

        \item $U$ is a stationary subset of $[H_\theta]^\omega$.
    \end{assume}
    Then in $V^{\col_U}$, for all limit $\eta<\alpha_g$, $\M_g(\eta)=\bigcup_{\xi<\eta}\M_g(\xi)$.
\begin{proof0}
    \item Assume otherwise.
    Then there exist $p$, $\alpha$, $\eta$, $x$, $M$ such that
    \begin{parts}
        \item $p\Vdash \alpha_g=\alpha$,

        \item $\eta$ is a limit ordinal strictly less than $\alpha$,

        \item $x\in H_\theta$,

        \item $p\Vdash x\in \M_g(\eta)-\bigcup_{\xi<\eta}\M_g(\xi)$,

        \item $p\Vdash \M_g(\eta)=M$,

        \item $M\in\M_p$.
    \end{parts}

    \item There exist $q\leq p$ and $N\in\M_q\cap M$ such that $x\in d_q(N)$.

    \item Let $r\leq q$ and $\xi<\eta$ be such that $r\Vdash \M_g(\xi)=N$.
    We have
    $$r\Vdash x\in d_r(N)\subseteq \M_g(\xi+1),$$
    which is a contradiction.
\end{proof0}

\item\proposition
    Suppose that 
    \begin{assume}
        \item $\theta>\omega$ is regular,

        \item $U$ is stationary in $[H_\theta]^\omega$.
    \end{assume}
    Then $\col_U$ is strongly semiproper.
\begin{proof0}
    \item Let us consider arbitrary
    \begin{parts}
        \item $\chi\gg\theta$,

        \item $M\prec (H_\chi,\in,U)$ countable,

        \item $p\in \col_U\cap M$.
    \end{parts}
    We want to find $q\leq p$ which is strongly semigeneric for $(M,\col_U)$.

    \item\textbf{Case I.} There exists $N\in U$ such that $M\cap H_\theta\subseteq N$ and $\delta_N=\delta_M$.
    \begin{proof0}
        \item Let $q:=(\M_p\cup\{N\},d_p\cup\{(N,\emptyset)\})$.
        We want to show that $q$ is strongly semigeneric for $(M,\col_U)$.

        \item Let $r\leq q$ be arbitrary and let us find some $r_M\in\col_U$ such that for $P:=\hull(M,r_M)$,
        \begin{parts}[ref=\arabic{pfenumii}$'$\alph{partsi}]
            \item $\delta_P=\delta_M$,

            \item\label{768'} for all $s\in\col_U\cap P$, if $s\leq r_M$, then $s\parallel r$.
        \end{parts}

        \item Let $r_M:=(\M_r\cap N,d_r\rest N)\in N$.
        We have that
        $$\hull(M,r_M)\cap H_\theta=\hull(M\cap H_\theta,r_M)\subseteq N.$$

        \item In particular, for $P:=\hull(M,r_M)$, we have $\delta_P=\delta_M$.

        \item Property \ref{768'}. is now routinely verified.
    \end{proof0}

    \item\textbf{Case II.} There does not exist $N\in U$ such that $M\cap H_\theta\subseteq N$ and $\delta_N=\delta_M$.
    \begin{proof0}
        \item In this case, we have that $q:=(\M_p\cup\{M\cap H_\theta\},d_p\cup\{(M\cap H_\theta,\emptyset)\}\in \col_U$ and $q\leq p$.
        We want to verify that $q$ is strongly generic for $(M,\col_U)$.

        \item Let $r\leq q$ be arbitrary and let us find some $r_M\in\col_U\cap M$ such that for all $s\in\col_U\cap M$, if $s\leq r_M$, then $s\parallel r$.

        \item It is routinely verified that $r_M:=(\M_r\cap M,d_r\rest M)$ is as required.
    \end{proof0}

    \item Having verified the two cases, we reach the conclusion.
\end{proof0}

\item\corollary
    Suppose that
    \begin{assume}
        \item $\theta>\omega$ is regular,

        \item $U$ is stationary in $[H_\theta]^\omega$,

        \item $g\leadsto V^{\col_U}$.
    \end{assume}
    Then $\otp(\M_g,\in)=\alpha_g=\omega_1$.
\begin{proof}
    This follows from the facts that $\alpha_g\leq\omega_1^V=\omega_1$ and $H_\theta^V=\cup\M_g$.
\end{proof}

\item\corollary
    Suppose that
    \begin{assume}
        \item $\theta>\omega$ is regular,

        \item $U$ is stationary in $[H_\theta]^\omega$.
    \end{assume}
    Then in $V^{\col_U}$, $|H_\theta^V|=\omega_1$.\qed

\item\lemma
\label{716}
    Suppose that
    \begin{assume}
        \item $\theta>\omega$ is regular,

        \item $U$ is stationary in $[H_\theta]^\omega$,

        \item $\chi\gg H_\theta$,

        \item $M\prec (H_\chi,\in,U)$ is countable,

        \item $p\in\col_U$ is such that $M\cap H_\theta\in\M_p$.
    \end{assume}
    Then $p$ is strongly generic for $(M,\col_U)$.
\begin{proof0}
    \item Let $q\leq p$ be arbitrary.
    We want to find $q_M\in\col_U\cap M$ such that for all $r\in\col_U\cap M$, if $r\leq q_M$, then $r\parallel q$.

    \item Let $q_M:=(\M_q\cap M,d_q\rest M)$.
    It is easily seen that $q_M$ is as required.
\end{proof0}

\item\proposition
    Suppose that
    \begin{assume}
        \item $\theta>\omega$ is regular,

        \item $U$ is stationary in $[H_\theta]^\omega$.
    \end{assume}
    Then in $V^{\col_U}$, $U$ is stationary in $[H_\theta^V]^\omega$.
\begin{proof0}
    \item Assume otherwise.
    Then there exist $C\in V^{\col_U}$ and $p\in\col_U$ such that 
    \begin{parts}
        \item $\Vdash_{\col_U}``C$ is a club in $[H_\theta^V]^\omega"$,

        \item $p\Vdash_{\col_U}C\cap U=\emptyset$.
    \end{parts}

    \item Let $\chi\gg H_\theta,C$.
    Then there exists countable $M\prec (H_\chi,\in,\theta,C,p)$ such that $M\cap H_\theta\in U$.

    \item Let $N:=M\cap H_\theta\in U$ and let $q:=(\M_p\cup\{N\},d_p\cup\{(N,\emptyset)\})\leq p$ in $\col_U$.
    By Lemma \ref{716}, we have that $q$ is generic for $(M,\col_U)$.

    \item Let $g\leadsto V^{\col_U}$ contain $q$.
    We have that $C\in M[g]$, so $M[g]\cap H_\theta^{V}\in C$.

    \item Since $q$ is generic for $(M,\col_U)$ in $V$, we have that $M[g]\cap H_\theta^V=N$.

    \item Thus, $N\in C\cap U$, which is a contradiction.
\end{proof0}
\end{enumerate}

\subsection{Sealed Predense Collections}
\begin{enumerate}

\item\definition
    Let $D\subseteq\P_\ns$.
    We say that \textit{$D$ is sealed} if there exists a mapping $f:\omega_1\to D$ and a club $C$ in $\omega_1$ such that $C\subseteq\nabla f$.

\item\proposition
    If $D\subseteq\P_\ns$ is sealed, then $D$ is predense in $\P_\ns$.\qed

\item\proposition
    The following are equivalent.
    \begin{parts}
        \item $\ns$ is saturated.

        \item Every predense subset of $\P_\ns$ is sealed.\qed 
    \end{parts}

\item\proposition
    Suppose that $D\subseteq\P_\ns$ is sealed and that $\P$ is a stationary set preserving poset.
    Then $V^\P\models``D$ is sealed in $\P_\ns"$.

\item\definition
    Suppose that
    \begin{assume}
        \item $D$ is predense in $\P_\ns$,

        \item $\theta\geq\omega_2$ is regular,

        \item $M\prec H_\theta$ is countable.
    \end{assume}
    Then we say that \textit{$M$ captures $D$} iff there exists $S\in D\cap M$ such that $\delta_M\in S$.

\item\lemma
    Suppose that
    \begin{assume}
        \item $D$ is predense in $\P_\ns$,

        \item $\theta\geq\omega_2$ is regular,

        \item $U:=\{M\prec H_\theta : |M|=\omega,``M\mbox{ captures } D"\}$.
    \end{assume}
    Then $U$ is projectively stationary.
\begin{proof0}
    \item Let $S\subseteq\omega_1$ be stationary and let $F:H_\theta^{<\omega}\to H_\theta$.
    There exists $T\in D$ such that $T\cap S$ is stationary.

    \item There exists $M\prec H_\theta$ such that $T\in M$, $F[M^{\omega}]\subseteq M$, and $\delta_M\in T\cap S$.

    \item Hence, $M\in U$, $\delta_M\in S$, and $M$ is closed for $F$, as required.
\end{proof0}

\item\proposition
\label{getting sealing}
    Suppose that
    \begin{assume}
        \item $D$ is a predense set in $\P_\ns$,

        \item $U:=\{M\prec H_{\omega_2} : |M|=\omega,``M\mbox{ captures }D"\}$,

        \item $g$ is $V$-generic for $\col_U$,

        \item $U^*:=\{\xi<\omega_1 : \M_g(\xi)\in U\}$,

        \item $h$ is generic over $V[g]$ such that $\omega_1^{V[g][h]}=\omega_1$,

        \item in $V[g][h]$, $U^*$ contains a club.
    \end{assume}
    Then in $V[g][h]$, $D$ is sealed.

\begin{proof0}
    \item For all $\xi<\omega_1$, we enumerate $D\cap \M_g(\xi)$ as $(S^\xi_n : n<\omega)$.

    \item We define $f:\omega_1\to D$ as follows: for all $\alpha<\omega_1$, for the unique $\xi<\omega_1$ and $n<\omega$ such that $\alpha=\omega\xi+n$, we let $f(\alpha):=S^\xi_n$.

    \item Let $C$ be a club contained in $U^*$.
    By Proposition \ref{continuity''}, the set $E:=\{\delta(\M_g(\alpha)) : \alpha\in C, \omega\alpha=\alpha\}$ is a club.

    \item We will be done if we show that $E\subseteq\nabla f$.
    To that end, let $\delta\in E$ be arbitrary and let us find $\gamma<\delta$ such that $\delta\in f(\gamma)$.

    \item There exists $\alpha\in C$ such that $\omega\alpha=\alpha$ and $\delta=\delta(\M_g(\alpha))$.
    
    \item Since $\alpha\in U^*$, we have that $\M_g(\alpha)\in U$, so there exists $S\in D\cap \M_g(\alpha)$ such that $\delta\in S$.

    \item Again by Proposition \ref{continuity''}, there exist $\xi<\alpha$ and $n<\omega$ such that $S=S^\xi_n$.

    \item Let $\gamma:=\alpha\xi+n<\alpha$.
    Then $\delta\in S^\xi_n=f(\gamma)$, as required.
\end{proof0}
\end{enumerate}

\subsection{Sealing Iteration}
\begin{enumerate}
\item\definition
    Let $\delta>\omega$ be a cardinal and let $A\subseteq V_\delta$.
    Then $\kappa<\delta$ is \textit{$A$-reflecting below $\delta$} if for all multiplicatively closed $\lambda\in(\kappa,\delta)$, there exists a strong $(\kappa,\lambda)$-extender $E$ satisfying $i_E(A\cap V_\kappa)\cap V_\lambda=A\cap V_\lambda$.

\item\definition
    Let $\delta>\omega$ be a cardinal.
    A \textit{Woodin diamond for $\delta$} is a function $\mathbf{U}:\kappa\to V_\kappa$ such that 
    \begin{parts}
        \item for all $\xi<\delta$, $\mathbf{U}(\xi)\subseteq V_\xi$,

        \item for all $A\subseteq V_\delta$, the set
        $$\{\kappa<\delta : A\cap V_\kappa=\uu(\kappa)\mbox{ and }\kappa\mbox{ is }A\mbox{-reflecting below }\delta\}$$
        is stationary in $\delta$.
    \end{parts}

\item\proposition
    Suppose that there exists a Woodin diamond for $\delta$.
    Then $\delta$ is a Woodin cardinal.\qed

\item\proposition
    Suppose that $\delta$ is a Woodin cardinal.
    Then in $V^{\col(\delta,\delta)}$, there exists a Woodin diamond for $\delta$.
\begin{proof}
    Cf. \cite[Lemma 1.3]{schindler2016nsomega1}
\end{proof}

\item\definition
    Let $\mathbf{U}$ be a Woodin diamond for $\delta$.
    \textit{The blueprint for the $\ns$-saturating iteration given by $\uu$} is the pair $((V_\delta,\in,\uu),\qq)$ where $\qq:\delta\times V_\delta\to V_\delta$ is defined over $(V_\delta,\in,\uu)$ as follows: for all $\kappa<\delta$ and for all $\P\in V_\delta$, if it is the case that
    \begin{parts}
        \item $\kappa$ is inaccessible,
        
        \item $\P\subseteq V_\kappa$ is a poset,
    
        \item\label{766} $\uu(\kappa)=\P\oplus (\dot S_i : i<\kappa)$,
    
        \item\label{768} $V^\P\models``\{\dot S_i : i<
        \kappa\}$ is predense in $\P_\ns"$,

        \item $\dot U$ is the canonical name such that 
        $$V^\P\models\dot U=\{M\prec H_{\omega_2} : |M|=\omega, ``M\mbox{ captures }\{\dot S_i : i<\kappa\}"\},$$
    \end{parts}
    then we set $\Q_\kappa$ to be the canonical $\P$-name for $\col_{\dot U}$, while otherwise we set $\Q_\kappa$ to be the canonical name for the trivial poset.

\item\proposition
    Suppose that
    \begin{assume}
        \item $\uu$ is a Woodin diamond for $\delta$,

        \item $(\V,\qq)$ is the blueprint for the $\ns$-saturating iteration given by $\uu$,

        \item $(\P_\alpha : \alpha\in\e^*)$ is the semi-proper forcing iteration given by the blueprint $(\V,\qq)$,

        \item $\P:=\bigcup_{\alpha\in\e^*}\P_\alpha$.
    \end{assume}
    Then it holds that
    \begin{parts}
        \item\label{1038} $\P$ is semi-proper,

        \item $\P$ is $\delta$-c.c,

        \item\label{1042} $\omega_2^{V^\P}=\delta$,

        \item in $V^\P$, $\ns$ is saturated.
    \end{parts}

\begin{proof0}
    \item Parts \ref{1038}-\ref{1042}. are general facts about the iteration.
    Let us verify that in $V^\P$, $\ns$ is saturated.

    \item Let $(\dot S_i : i<\delta)$ be a name for a predense collection in $\P_\ns$ in $V^\P$.
    We want to show that the corresponding predense collection in $V^\P$ is sealed.

    \item The set
    $$\e_0:=\left\{\alpha<\delta : V_\alpha\prec (V_\delta,\in,\uu,\P,(\dot S_i : i<\delta))\right\}$$
    is a club in $\delta$.
    This means that there exists $\kappa\in\e_0$ such that 
    \begin{parts}
        \item $\uu(\kappa)=\P_{<\kappa}\oplus (\dot S_i : i<\kappa)$,

        \item $\kappa$ is $\P\oplus (S_i : i<\delta)$-reflecting below $\delta$.
    \end{parts}
    We will show that $(\dot S_i : i<\kappa)$ is sealed inside $V^{\P_\kappa}$.

    \item Let $\lambda,E,j,\Mod$ be such that
    \begin{parts}
        \item $\lambda\in\e_0$ is inaccessible $>\kappa$,


        \item $E$ is a strong $(\kappa,\lambda)$-extender,

        \item $j:=i_E:V\prec\ult(V,E)=:\Mod$,

        \item $j(\P_{<\kappa})\cap V_\lambda=\P_{<\lambda}$,

        \item $j((\dot S_i : i<\kappa))\rest\lambda=(\dot S_i : i<\lambda)$.
    \end{parts}

    \item Let $G_\kappa$ be $V$-generic for $\P_\kappa$ and let us work in $V[G_\kappa]$.
    The family $(S_i : i<\kappa)$ is predense in $\P_\ns$, so $\Q_\kappa=\col_U$, where
    $$U:=\{M\prec H_{\omega_2} : |M|=\omega,``M\mbox{ captures }\{S_i : i<\kappa\}"\}.$$

    \item Let $g$ be $V[G_\kappa]$-generic for $\col_U$, let $G_{\kappa+1}$ be the $V$-generic for $\P_{\kappa+1}$ corresponding to $G_\kappa*g$, and let us work in $V[G_{\kappa+1}]$.
    We denote by $U^*$ the set
    $$\{\xi<\omega_1 : \M_g(\xi)\in U\}.$$

    \item Let $G^*$ be $V$-generic for $j(\P_\kappa)$ extending $G_{\kappa+1}$ and let $k:V[G_\kappa]\prec \Mod [G^*]$ be such that $k\rest V=j$ and $k(G_\kappa)=G^*$.

    \item\claim In $\Mod [G^*]$, $\omega_1-U^*$ is not stationary.
    \begin{proof0}
        \item Let us assume otherwise.
        Then there exists $i_0<j(\kappa)$ such that the set
        $$k((\dot S_i : i<\kappa))(i_0)\cap (\omega_1-U^*)$$
        is stationary in $\omega_1$.

        \item In $\Mod$, let $(\dot T_i : i<j(\kappa))=k((\dot S_i : i<\kappa))$ and let $(M_\xi : \xi<\omega_1)$ be a continuous $\in$-chain of countable elementary submodels of
        $$(\Mod\rest j(\lambda),\in,\kappa,\lambda,\dot T_i).$$

        \item Back in $\Mod [G^*]$, the set
        $$\{\delta(M_\xi [G^*]) : \xi<\omega_1\}$$
        is a club in $\omega_1$, so there exists $\xi_0<\omega_1$ such that
        $$\alpha:=\delta(M_{\xi_0}[G^*])\in T_{i_0}\cap(\omega_1-U^*).$$

        \item Since $\alpha\not\in U^*$, we have that $\M_g(\alpha)\not\in U$.
        In other words, inside $\Mod [G_\kappa]$, there does not exists $N\in U$ such that $N\supseteq\M_g(\alpha)$ and $\delta_N=\delta(\M_g(\alpha))$.

        \item Since $\M_g(\alpha)$ is a countable subset of $H_{\omega_2}$ inside $\Mod [G_\kappa]$ and since $\omega_2^{\Mod [G_\kappa]}=\kappa$, we have that $k(\M_g(\alpha))=\M_g(\alpha)$.

        \item Consequently, in $\Mod [G^*]$, there does not exist $N\in k(U)$ such that $N\supseteq \M_g(\alpha)$ and $\delta_N=\delta(\M_g(\alpha))$.

        \item Since $g\in M_{\xi_0}[G^*]$, we get that
        $$M_{\xi_0}[G^*]\cap H_{\omega_2}^{V[G_\kappa]}=\bigcup_{\xi<\alpha}\M_g(\xi)=\M_g(\alpha).$$
        Hence, $M_{\xi_0}[G^*]\supseteq\M_g(\alpha)$ and $\delta(M_{\xi_0}[G^*])=\delta(\M_g(\alpha))$.

        \item On the other hand, $M_{\xi_0}[G^*]$ captures $k(\{S_i : i<\kappa\})$, as witnessed by $T_{i_0}$.
        This means that $M_{\xi_0}[G^*]\in k(U)$.
        
        \item The previous three points together yield a contradiction.
    \end{proof0}

    \item Thus, in $\Mod [G^*]$, $U^*$ contains a club.
    
    \item By Proposition \ref{getting sealing} applied in $\Mod [G_\kappa]$, we have that $\{S_i : i<\kappa\}$ is sealed in $\Mod [G^*]$.

    \item This means that $k(\{S_i : i<\kappa\})$ is sealed in $\Mod [G^*]$, which by elementarity of $k$ yields that $\{S_i : i<\kappa\}$ is sealed in $V[G_\kappa]$, as required.
\end{proof0}

\item\corollary
    Suppose that $\delta$ is a Woodin cardinal.
    Then there exists a poset $\P$ such that
    \begin{parts}
        \item $\P$ is semi-proper,

        \item $\P$ is $\delta$-c.c,

        \item $\omega_2^{V^\P}=\delta$,

        \item in $V^\P$, $\ns$ is saturated.\qed
    \end{parts}
\end{enumerate}
\newpage

\nocite{velickovic2021iteration}
\bibliographystyle{alpha}
\bibliography{lit.bib}
\end{document}